\documentclass[11pt]{amsart}

\usepackage{amsfonts}
\usepackage{mathrsfs}
\usepackage{amscd}
\usepackage{amsmath}
\usepackage{amssymb}
\usepackage{graphicx}
\usepackage{latexsym}
\usepackage{lscape}
\usepackage{xypic}
\usepackage{comment}
\usepackage{amscd}
\usepackage{wasysym}  
\usepackage{stmaryrd}  
\usepackage{tikz} 
\usetikzlibrary{shapes.geometric, calc, arrows.meta}
\usepackage{tikz-cd}
\usetikzlibrary{matrix,arrows}
\usepackage{appendix}
\usepackage{geometry}
\usepackage{needspace}
\usepackage[nice]{nicefrac}
\usepackage{dsfont}
\usepackage[]{overpic}

\usepackage[pagebackref,hyperindex,linktocpage=true]{hyperref}
\hypersetup{
    colorlinks,
    linkcolor={red!50!black},
    citecolor={blue!50!black},
    urlcolor={blue!80!black}
}

\usepackage{color}

\newtheorem{thm}{Theorem}[section]

\newtheorem{prop}[thm]{Proposition}
\newtheorem{lem}[thm]{Lemma}

\newtheorem{lem-def}[thm]{Lemma-Definition}
\newtheorem{cor}[thm]{Corollary}

\newtheorem{thmx}{Theorem}

\theoremstyle{definition}

\newtheorem{ex}[thm]{Example}

\newtheorem{rmk}[thm]{Remark}
\newtheorem{dfn}[thm]{Definition}
\newtheorem*{prop*}{Proposition}

\numberwithin{equation}{section}

\newcommand{\nc}{\newcommand}

\nc{\on}{\operatorname}
\nc{\fraka}{{\mathfrak a}} \nc{\bba}{{\mathbf a}}
\nc{\frakb}{{\mathfrak b}}
\nc{\frakc}{{\mathfrak c}}
\nc{\frakd}{{\mathfrak d}}
\nc{\frake}{{\mathfrak e}}
\nc{\frakf}{{\mathfrak f}}
\nc{\frakg}{{\mathfrak g}}
\nc{\frakh}{{\mathfrak h}}
\nc{\fraki}{{\mathfrak i}}
\nc{\frakj}{{\mathfrak j}}
\nc{\frakk}{{\mathfrak k}}
\nc{\frakl}{{\mathfrak l}}
\nc{\frakm}{{\mathfrak m}}
\nc{\frakn}{{\mathfrak n}}
\nc{\frako}{{\mathfrak o}}
\nc{\frakp}{{\mathfrak p}}
\nc{\frakq}{{\mathfrak q}}
\nc{\frakr}{{\mathfrak r}}
\nc{\fraks}{{\mathfrak s}}
\nc{\frakt}{{\mathfrak t}}
\nc{\fraku}{{\mathfrak u}}
\nc{\frakv}{{\mathfrak v}}
\nc{\frakw}{{\mathfrak w}}
\nc{\frakx}{{\mathfrak x}}
\nc{\fraky}{{\mathfrak y}}
\nc{\frakz}{{\mathfrak z}}
\nc{\frakA}{{\mathfrak A}}
\nc{\frakB}{{\mathfrak B}}
\nc{\frakC}{{\mathfrak C}}
\nc{\frakD}{{\mathfrak D}}
\nc{\frakE}{{\mathfrak E}}
\nc{\frakF}{{\mathfrak F}}
\nc{\frakG}{{\mathfrak G}}
\nc{\frakH}{{\mathfrak H}}
\nc{\frakI}{{\mathfrak I}}
\nc{\frakJ}{{\mathfrak J}}
\nc{\frakK}{{\mathfrak K}}
\nc{\frakL}{{\mathfrak L}}
\nc{\frakM}{{\mathfrak M}}
\nc{\frakN}{{\mathfrak N}}
\nc{\frakO}{{\mathfrak O}}
\nc{\frakP}{{\mathfrak P}}
\nc{\frakQ}{{\mathfrak Q}}
\nc{\frakR}{{\mathfrak R}}
\nc{\frakS}{{\mathfrak S}}
\nc{\frakT}{{\mathfrak T}}
\nc{\frakU}{{\mathfrak U}}
\nc{\frakV}{{\mathfrak V}}
\nc{\frakW}{{\mathfrak W}}
\nc{\frakX}{{\mathfrak X}}
\nc{\frakY}{{\mathfrak Y}}
\nc{\frakZ}{{\mathfrak Z}}
\nc{\bbA}{{\mathbb A}}
\nc{\bbB}{{\mathbb B}}
\nc{\bbC}{{\mathbb C}}
\nc{\bbD}{{\mathbb D}}
\nc{\bbE}{{\mathbb E}}
\nc{\bbF}{{\mathbb F}} \nc{\bbf}{{\mathbf f}}
\nc{\bbG}{{\mathbb G}}
\nc{\bbH}{{\mathbb H}}
\nc{\bbI}{{\mathbb I}}
\nc{\bbJ}{{\mathbb J}}
\nc{\bbK}{{\mathbb K}}
\nc{\bbL}{{\mathbb L}}
\nc{\bbM}{{\mathbb M}}
\nc{\bbN}{{\mathbb N}}
\nc{\bbO}{{\mathbb O}}
\nc{\bbP}{{\mathbb P}}
\nc{\bbQ}{{\mathbb Q}}
\nc{\bbR}{{\mathbb R}}
\nc{\bbS}{{\mathbb S}}
\nc{\bbT}{{\mathbb T}}
\nc{\bbU}{{\mathbb U}}
\nc{\bbV}{{\mathbb V}}
\nc{\bbW}{{\mathbb W}}
\nc{\bbX}{{\mathbb X}}
\nc{\bbY}{{\mathbb Y}}
\nc{\bbZ}{{\mathbb Z}}
\nc{\calA}{{\mathcal A}}
\nc{\calB}{{\mathcal B}}
\nc{\calC}{{\mathcal C}}
\nc{\calD}{{\mathcal D}}
\nc{\calE}{{\mathcal E}}
\nc{\calF}{{\mathcal F}}
\nc{\calG}{{\mathcal G}}
\nc{\calH}{{\mathcal H}}
\nc{\calI}{{\mathcal I}}
\nc{\calJ}{{\mathcal J}}
\nc{\calK}{{\mathcal K}}
\nc{\calL}{{\mathcal L}}
\nc{\calM}{{\mathcal M}}
\nc{\calN}{{\mathcal N}}
\nc{\calO}{{\mathcal O}}
\nc{\calP}{{\mathcal P}}
\nc{\calQ}{{\mathcal Q}}
\nc{\calR}{{\mathcal R}}
\nc{\calS}{{\mathcal S}}
\nc{\calT}{{\mathcal T}}
\nc{\calU}{{\mathcal U}}
\nc{\calV}{{\mathcal V}}
\nc{\calW}{{\mathcal W}}
\nc{\calX}{{\mathcal X}}
\nc{\calY}{{\mathcal Y}}
\nc{\calZ}{{\mathcal Z}}

\nc{\scrA}{{\mathscr A}}
\nc{\scrB}{{\mathscr B}}
\nc{\scrR}{{\mathscr R}}

\nc{\blambda}{{\bar{ \lambda}}}

\nc{\olO}{\bar{\calO}}

\nc{\al}{{\alpha}} 
\nc{\be}{{\beta}}
\nc{\ga}{{\gamma}} \nc{\Ga}{{\Gamma}}
 \nc{\hGa}{\widehat{\Gamma}}
\nc{\ve}{{\varepsilon}} 
\nc{\la}{{\lambda}} \nc{\La}{{\Lambda}}
\nc{\om}{\omega} \nc{\Om}{\Omega} 
\nc{\sig}{{\sigma}} \nc{\Sig}{{\Sigma}}
\nc{\de}{\delta} \nc{\De}{\Delta}

\nc{\tnb}{\psi_{\rm tame}}
\nc{\oM}{\overline{{M}}}
\nc{\op}{{\on{op}}}
\nc{\ad}{{\on{ad}}}
\nc{\alg}{{\on{alg}}}
\nc{\Ad}{{\on{Ad}}}
\nc{\Adm}{{\on{Adm}}} \nc{\aff}{{\on{af}}}
\nc{\Aut}{{\on{Aut}}}
\nc{\Bun}{{\on{Bun}}}
\nc{\cha}{{\on{char}}}
\nc{\der}{{\on{der}}}
\nc{\Der}{{\on{Der}}}
\nc{\diag}{{\on{diag}}}
\nc{\End}{{\on{End}}}
\nc{\Fl}{{\calF\ell}}
\nc{\Tr}{{\on{Transp}}}
\nc{\TR}{{\calT\calR}}
\nc{\Gal}{{\on{Gal}}}
\nc{\Gr}{{\on{Gr}}}
\nc{\rH}{{\on{H}}}
\nc{\Hom}{{\on{Hom}}}
\nc{\IC}{{\on{IC}}}
\nc{\id}{{\on{id}}}
\nc{\Id}{{\on{Id}}}
\nc{\ind}{{\on{ind}}}
\nc{\Ind}{{\on{Ind}}}
\nc{\Lie}{{\on{Lie}}}
\nc{\Pic}{{\on{Pic}}}
\nc{\pr}{{\on{pr}}}
\nc{\Res}{{\on{Res}}}
\nc{\res}{{\on{res}}} \nc{\Sat}{{\on{Sat}}}
\nc{\s}{{\on{sc}}}
\nc{\drv}{{\on{der}}}
\nc{\sgn}{{\on{sgn}}}
\nc{\Spec}{{\on{Spec}}}\nc{\Spf}{\on{Spf}} 
\nc{\Sph}{\on{Sph}}
\nc{\St}{{\on{St}}}
\nc{\tr}{{\on{tr}}}
\nc{\Mod}{{\mathrm{-Mod}}}
\nc{\Hilb}{{\on{Hilb}}} 
\nc{\Ext}{{\on{Ext}}} 
\nc{\vs}{{\on{Vec}}}
\nc{\ev}{{\on{ev}}}
\nc{\nO}{{\breve{\calO}}}
\nc{\tS}{{\widetilde{S}}}
\nc{\spe}{{\on{sp}}}
\nc{\loc}{{\on{loc}}}

\nc{\nscrR}{{\mathscr{R}^{\on{nr}}}}

\nc{\GL}{{\on{GL}}}
\nc{\U}{{\on{U}}}
\nc{\Gl}{\on{Gl}} 
\nc{\GSp}{{\on{GSp}}}
\nc{\gl}{{\frakg\frakl}}
\nc{\SL}{{\on{SL}}} 
\nc{\SU}{{\on{SU}}} 
\nc{\SO}{{\on{SO}}}
\nc{\PGL}{{\on{PGL}}}

\nc{\Conv}{{\on{Conv}}}
\nc{\Rep}{{\on{Rep}}}
\nc{\Dom}{{\on{Dom}}}
\nc{\red}{{\on{red}}}
\nc{\act}{{\on{act}}}
\nc{\nr}{{\on{nr}}}
\nc{\ctf}{{\on{ctf}}}

\nc{\str}{{\on{-}}} 
\nc{\os}{{\bar{s}}}
\nc{\oeta}{{\bar{\eta}}}

\nc{\hookto}{\hookrightarrow}
\nc{\longto}{\longrightarrow}
\nc{\leftto}{\leftarrow}
\nc{\onto}{\twoheadrightarrow}
\nc{\lonto}{\twoheadleftarrow}

\nc{\uG}{{\underline{G}}}
\nc{\uA}{{\underline{A}}}
\nc{\uS}{{\underline{S}}}
\nc{\uT}{{\underline{T}}}
\nc{\uM}{{\underline{M}}}
\nc{\uP}{{\underline{P}}}
\nc{\uB}{{\underline{B}}}
\nc{\uN}{{\underline{N}}}

\nc{\ucG}{{\underline{\calG}}}
\nc{\ucA}{{\underline{\calA}}}
\nc{\ucS}{{\underline{\calS}}}
\nc{\ucT}{{\underline{\calT}}}
\nc{\ucM}{{\underline{\calM}}}
\nc{\ucP}{{\underline{\calP}}}
\nc{\ucN}{{\underline{\calN}}}

\nc{\bF}{{\breve{F}}}

\nc{\oFl}{{\overline{\Fl}}} 
\nc{\bU}{{\overline{U}}}
\nc{\tGr}{{\widetilde{\Gr}}}
\nc{\cGr}{\calG r}
\nc{\oGr}{\overline{\on{Gr}}} 
\nc{\ocGr}{\overline{\calG r}}
\nc{\co}{{\colon}}
\nc{\sch}[1]{(Sch/{#1})}
\nc{\HypLoc}[1]{HypLoc({#1})}

\nc{\ohtimes}{\stackrel{!}{\otimes}}
\nc{\boxtilde}{\widetilde{\boxtimes}}
\nc{\vstar}{{\varhexstar}}

\nc{\Div}{\on{Div}}

\nc{\bslash}{\backslash}
\nc{\algQl}{{\bar{\bbQ}_\ell}}
\nc{\sF}{{\bar{F}}}
\nc{\nF}{{\breve{F}}}
\nc{\nW}{{W^{\on{nr}}}}
\nc{\sk}{{\bar{k}}}
\nc{\cont}{\on{c}}
\nc{\Supp}{\on{Supp}}
\nc{\blt}{\bullet}  
\nc{\dom}{\on{dom}}
\nc{\scon}{{\on{sc}}} 
\nc{\Affine}{\on{Aff}} 
\nc{\nscrA}{\mathscr{A}^{\on{nr}}} 
\nc{\nfraka}{{\bbf^{\on{nr}}}}
\nc{\ran}{{\rangle}}
\nc{\lan}{{\langle}}
\nc{\bk}{{\bar{k}}}
\nc{\tF}{{\widetilde{F}}}
\nc{\sS}{{\bar{S}}}
\nc{\LG}{{^\text{L}\hspace{-0.04cm}G}}
\nc{\LL}{{^\text{L}\hspace{-0.07cm}L}}

\nc{\pot}[1]{ [\hspace{-0,5mm}[ {#1} ]\hspace{-0,5mm}] }
\nc{\rpot}[1]{ (\hspace{-0,7mm}( {#1} )\hspace{-0,7mm}) }

\nc{\defined}{\hspace{0.1cm}\stackrel{\text{\tiny \rm def}}{=}\hspace{0.1cm}}

\begin{document}

\title[Alcove walks, crystal structure and a Levi branching rule]{Alcove walks, crystal structure and a Levi branching rule}
\author[W.~Wang]{by Wenzhuo Wang}

\address{Department of Mathematics, University of Maryland, College Park, MD 20742-4015, USA}
\email{wzwang@umd.edu}

\begin{abstract}
In the Borel case \(P=B\), we endow Haines' set of maximal-dimensional positively folded alcove walks
with a crystal structure. Furthermore, we construct a bijection
between this alcove-walk model and the Gaussent--Littelmann model of
LS galleries, and we prove that the bijection is compatible with the
root operators. Thus the crystal structure on the alcove-walk model is
obtained by transport from the Gaussent--Littelmann crystal. This answers
a question raised in \cite{Hai24}. We also generalize the stratification result in \cite{GL05} to the parabolic setting, give a proof of the Levi branching rule for LS galleries, and compare it with the general maximal-dimensional positively folded alcove walks in \cite{Hai24}.
\end{abstract}

\maketitle

\setcounter{tocdepth}{1}
\tableofcontents
\setcounter{section}{0}

\thispagestyle{empty}

\section{Introduction}

Let $k$ be a field, let $F=k((t))$, and let $\calO=k[[t]]$. Let $G$ be a connected split reductive group over $k$, and fix a maximal torus and a Borel subgroup $T\subset B\subset G$. Write $K=G(\calO)$. Let $\widehat G$ be the complex dual group, with dual Borel pair $\widehat T\subset\widehat B$. For every dominant cocharacter $\mu\in X_{*}(T)^+$, let $V_{\mu}^{\widehat G}$ denote the irreducible representation of $\widehat G$ of highest weight $\mu$.

Several combinatorial models describe the weight multiplicities of $V_{\mu}^{\widehat G}$. Littelmann's LS paths of shape $\mu$ form the path model \cite{Li95}. Gaussent and Littelmann introduced the LS-gallery model, which replaces paths by chains of adjacent facets in a fixed apartment of the Bruhat--Tits building $\calB(G,F)$ \cite{GL05}. Braverman and Gaitsgory studied the crystal structure on the set of top-dimensional irreducible components of the intersection
\begin{equation*}
LUt^{\lambda}K/K\cap Kt^{\mu}K/K.
\end{equation*}
These irreducible components are the MV cycles of coweight $(\lambda,\mu)$, and their number is the dimension of the $\lambda$-weight space of $V_{\mu}^{\widehat G}$ \cite{BG01}. Finally, for a standard parabolic subgroup $P=MN$, Haines used retractions and alcove walks in $\calB(G,F)$ to construct a combinatorial model for the top-dimensional irreducible components of the parabolic MV intersection
\begin{equation*}
L^+MLN t^{\lambda}K/K\cap Kt^{\mu}K/K
\end{equation*}
\cite{Hai24}.

An alcove walk moves through adjacent alcoves in a fixed apartment;
a fold means staying in the same alcove at a prescribed panel.
Haines' walks start at the base alcove $\Omega_0$.
Their types come from minimal walks from $\Omega_0$ to the nearest
alcoves whose closures contain vertices in $W_0(-\mu)$, where $W_0$ is the
finite Weyl group of $G$.
After folding, the terminal alcove is required to contain a vertex
of $-W_M\lambda$ in its closure, where $W_M$ is the Weyl group
of $M$. In the Borel case this vertex is $-\lambda$; see
\cite[Definition~5.8]{Hai24}.

In the Borel case $P=B$, Haines' set
$\mathcal M_\mu^{\De_U}(\lambda)$ of maximal-dimensional
positively folded alcove walks parametrizes the MV cycles
of coweight $(\lambda,\mu)$ \cite[Theorem~5.15]{Hai24}.
Our first goal is to equip this model with explicit raising
and lowering operators. The starting point is a geometric
map from alcove walks to LS galleries: the skeleton map.

Put $\mu^*=-w_0\mu$, where $w_0$ is the longest element of $W_0$.
Let $W_{\mu^*}$ be the stabilizer of $\mu^*$ in $W_0$ and let $W^{\mu^*}$ be the
minimal-length representatives of $W_0/W_{\mu^*}$.
Using the reduced expressions chosen in Section~4, a walk
$\delta'\in\mathcal M_\mu^{\De_U}(\lambda)$ has type
$w(t_{\mu^*})_0$ for some $w\in W^{\mu^*}$.
Lemma~\ref{lem:initial-minimal-segment} shows that its first
$\ell(w)$ steps form a minimal walk from $\Omega_0$ to $w\Omega_0$,
entirely among alcoves having $0$ in their closures.
Remove precisely these steps, retaining $w\Omega_0$ as the first
alcove, and mark $0$ and $-\lambda$ as the initial and terminal
vertices. The resulting vertex-to-vertex walk $\delta$ has the
fixed type $\de_\mu$. Thus Haines' model fixes the initial alcove
and allows the type to vary, whereas our model fixes the type
and allows the initial alcove to vary.
The skeleton $\operatorname{sk}(\delta)$ is obtained by retaining the
prescribed faces of $\delta$; see
Definition~\ref{def:skeleton}.
Write $\Ga_{LS}^+(\lambda,\ga_\mu)$ for the LS galleries
of the fixed type $\ga_\mu$ ending at $-\lambda$.
Our first result shows that passing to the skeleton loses
no information on Haines' model.

\begin{thmx}\textup{(}Proposition~\ref{prop:skeleton-bijection}\textup{)}
The skeleton map
\[
\begin{aligned}
\mathcal M_\mu^{\De_U}(\lambda)
&\longrightarrow \Ga_{LS}^+(\lambda,\ga_\mu),\\
\delta'&\longmapsto\operatorname{sk}(\delta),
\end{aligned}
\]
is a bijection.
\end{thmx}

To define the raising and lowering operators directly on
alcove walks, we introduce the auxiliary models
$B_\mu^v(\lambda)$ of Definition~\ref{def:alcove-model}.
Here $v\in W_{\mu^*}$. The model $B_\mu^e(\lambda)$ consists
of the remainders $\delta$ just described. More generally,
$B_\mu^v(\lambda)$ allows the first alcove to be $wv\Omega_0$,
with $w\in W^{\mu^*}$, while keeping the same walk type.
The additional choices are needed when $\mu$ is singular,
since the reflections in the root-operator formulas can
change the initial alcove. Reflections and translations
of subwalks, together with the braid-move comparisons
between these models, give the following result.

\begin{thmx}\textup{(}Theorem~\ref{thm:alcove-crystal} and
Section~\ref{subsec:explicit-formulas}\textup{)}
Put $B_\mu^e=\bigsqcup_\lambda B_\mu^e(\lambda)$.
For each simple root $\al$ of $G$, there are
combinatorially defined raising and lowering operators
\[
e_\al,f_\al:B_\mu^e\longrightarrow
B_\mu^e\cup\{\mathbf0\}.
\]
With $\operatorname{wt}(\de)=\lambda$ for
$\de\in B_\mu^e(\lambda)$, these operators make $B_\mu^e$
a highest-weight crystal of highest weight $\mu$.
The skeleton bijection intertwines them with the
Gaussent--Littelmann operators on LS galleries.
\end{thmx}

Our second goal is a gallery model for Levi branching.
For a standard parabolic subgroup $P=MN$, we decompose the
parabolic MV intersection into strata indexed by
$\De_{I_P}$-positively folded galleries.
Here $\De_{I_P}$ is the base alcove translated sufficiently
far toward infinity in the positive $N$-directions; see Section~2.3.
The strata of dimension $\langle\rho,\lambda+\mu\rangle$
are indexed precisely by the LS galleries lying in the closed
$M$-antidominant Weyl chamber. This gives the following
extension of \cite{GL05}.

\begin{thmx}\textup{(}Propositions~\ref{prop:parabolic-strata},
\ref{prop:parabolic-retraction-decomposition}, \ref{prop:parabolic-component-count} and
Theorem~\ref{thm:levi-gallery-branching}\textup{)}
Let $\mu$ be $G$-dominant and let $\lambda$ be $M$-dominant.
Put
\[
Z_P(\lambda)=L^+MLNt^\lambda K/K\cap Kt^\mu K/K.
\]
Let $\Ga_P^+(\nu,\ga_\mu)$ be the set of
$\De_{I_P}$-positively folded galleries of type $\ga_\mu$
ending at $-\nu$. There is a finite decomposition
\[
Z_P(\lambda)
=\bigsqcup_{\nu\in W_M\lambda}
 \bigsqcup_{\ga\in\Ga_P^+(\nu,\ga_\mu)}S_P(\ga)
\]
into geometrically irreducible locally closed strata, with
$\dim S_P(\ga)=\dim_P(\ga)$.
Here $S_P(\ga)$ is the parabolic retraction stratum of
Section~7.1, and $\dim_P$ is the gallery dimension for the
$\De_{I_P}$-orientation.
Taking closures of the strata of dimension
$\langle\rho,\lambda+\mu\rangle$ gives a bijection
\[
\begin{aligned}
\bigl\{\ga\in\Ga_{LS}^+(\lambda,\ga_\mu)
 \bigm|\ga\subset\overline{\calC_M^-}\bigr\}
&\xrightarrow{\sim}\operatorname{Irred}_{\langle\rho,\lambda+\mu\rangle} Z_P(\lambda),\\
\ga&\longmapsto\overline{S_P(\ga)}.
\end{aligned}
\]
Thus these LS galleries parametrize the parabolic MV cycles
of maximal possible dimension. In particular, the Levi
branching multiplicity is
\[
[V_\mu^{\widehat G}:V_\lambda^{\widehat M}]
=\#\bigl\{\ga\in\Ga_{LS}^+(\lambda,\ga_\mu)
 \bigm|\ga\subset\overline{\calC_M^-}\bigr\}.
\]
\end{thmx}

The decomposition is induced by the parabolic retraction
$\rho_{I_P,\calA}$, applied face by face to minimal galleries
of the fixed type starting at $0$. On these galleries it
agrees with the retraction centered at $\De_{I_P}$.
The stratum $S_P(\ga)$ consists of the endpoints of those
minimal galleries $\widetilde\ga$ satisfying
$\rho_{I_P,\calA}(\widetilde\ga)=\ga$; its nonemptiness is
equivalent to $\ga$ being $\De_{I_P}$-positively folded.
Proposition~\ref{prop:parabolic-retraction-decomposition} gives
the decomposition, and Proposition~\ref{prop:parabolic-strata}
computes the dimensions of its strata.
Proposition~\ref{prop:parabolic-component-count} identifies the closures of the
$\langle\rho,\lambda+\mu\rangle$-dimensional strata with the
parabolic MV cycles. Theorem~\ref{thm:levi-gallery-branching} identifies their indexing galleries with the LS
galleries contained in $\overline{\calC_M^-}$.
The proof follows the geometric strategy of \cite{GL05}.

The same LS galleries give a direct construction of all the
maximal-dimensional $\De_{I_P}$-positively folded alcove walks
in Haines' set $\mathcal M_\mu^{\De_{I_P}}(\lambda)$
\cite[Definition~5.8 and equation~(5.5)]{Hai24}.
No new lifting procedure is needed: we use the Borel inverse
skeleton construction of Section~\ref{subsec:haines-inverse}.

\begin{thmx}\textup{(}Proposition~\ref{prop:parabolic-walk-criterion}\textup{)}
\label{thm:intro-haines-walks}
Let $\mu$ be $G$-dominant and let $\lambda$ be $M$-dominant.
For each $\ga\in\Ga_{LS}^+(\lambda,\ga_\mu)$ contained in
$\overline{\calC_M^-}$, let $\delta=L(\ga)\in B_\mu^e(\lambda)$
be its Borel alcove-walk lift. Write $w\Omega_0$ for the
initial alcove of $\delta$, where $w\in W^{\mu^*}$.
Prepend the fixed $\ell(w)$-step minimal walk from $\Omega_0$ to
$w\Omega_0$, and denote the resulting walk by $\delta'$.
This construction gives a bijection
\[
\bigl\{\ga\in\Ga_{LS}^+(\lambda,\ga_\mu)
 \bigm|\ga\subset\overline{\calC_M^-}\bigr\}
\xrightarrow{\sim}\mathcal M_\mu^{\De_{I_P}}(\lambda),
\qquad \ga\longmapsto\delta'.
\]
In particular, every walk in Haines' set is obtained exactly
once and ends at $-\lambda$.
\end{thmx}

Explicitly, if $\ga=e_{\al_N}^{GL}\cdots e_{\al_1}^{GL}(\ga_\mu)$,
then $L(\ga)=e_{\al_N}\cdots e_{\al_1}(\de_\mu)$,
using the raising operators of Theorem~B.
The result is independent of the chosen sequence of raising
operators; see Proposition~\ref{prop:inverse-by-raising}.

\medskip
\noindent\textbf{Outline.}
Section~2 reviews the necessary group-theoretic background and retractions of the building. Section~3 recalls galleries and alcove walks, defines the root operators in our convention, compares them with those of \cite{GL05}, and proves the combinatorial lemmas needed later. Section~4 studies retractions on Bott--Samelson varieties and the maps between the resulting stratifications. Sections~5 and~6 construct the bijection between $B_{\mu}^v$ and $\Ga_{LS}^+(\ga_{\mu})$ and prove compatibility with the crystal operators. Section~7 proves the Levi branching rule and compares the resulting LS-gallery and alcove-walk models.

\medskip
\noindent\textbf{Acknowledgments.}
I am deeply grateful to my advisor, Thomas J. Haines, for suggesting this project, for many helpful discussions, and for his careful reading and detailed comments on an earlier draft. The author was supported in part by the National Science Foundation under Grant No. DMS-2200873 during the preparation of this work.

\medskip
\noindent\textbf{AI usage disclosure.}
Generative AI was used only to edit and proofread the article. All mathematical content was developed and verified by the author, who takes full responsibility for the content of this article.

\section{Background}

Every scheme and ind-scheme in this article is equipped with its reduced structure.

\subsection{Group-theoretic background}
Let $G$ be a connected split reductive group over a field $k$. We fix a Borel pair $T\subset B$ in $G$, which determines the dual Borel pair $\widehat T\subset\widehat B$ in the dual group $\widehat G$. Let $\Phi=\Phi(G,T)$ and $\Phi^\vee$ be the associated sets of roots and coroots. Write
$\Phi^+\subset\Phi$ for the $B$-positive roots, put
$\Phi^-=-\Phi^+$, and set
$\rho=\frac12\sum_{\alpha\in\Phi^+}\alpha$.
For a standard Levi subgroup $M$, write
$\Phi_M=\Phi(M,T)$ and $\Phi_M^+=\Phi_M\cap\Phi^+$.

Let $\calO=k[[t]]$ be the ring of formal power series and let $F=k((t))$ be its fraction field. The Iwahori--Weyl group
\begin{equation*}
W=N_G(T)(F)/T(\calO)
\end{equation*}
fits into the exact sequence
\begin{equation*}
1\longrightarrow T(F)/T(\calO)\longrightarrow W
\longrightarrow N_G(T)(F)/T(F)\longrightarrow 1.
\end{equation*}
The first term is naturally identified with $X_*(T)$, and the last quotient is the finite Weyl group $W_0$. The normalizer $N_G(T)(k)$ identifies $W_0$ with a subgroup of $W$; we use this splitting throughout. Thus we identify $W$ with $W_0\ltimes X_*(T)$. For $\lambda\in X_*(T)$, let $t^\lambda\in W$ be the image of $t$ under $\lambda:F^\times\to T(F)$. Let $\calA=\calA_T$ be the apartment associated with $T$ in the enlarged Bruhat--Tits building $\calB(G,F)$. We use the Tits convention: $t^\lambda$ acts on $\calA_T$ by translation by $-\lambda$. We also fix representatives in $N_G(T)(\calO)$ of the elements of $W_0$.

For every root $\al\in\Phi(G,T)$, let $U_\al$ be the corresponding root subgroup. Choose compatible root-group parametrizations $u_\al:\bbG_a\to U_\al$ extending the fixed pinning. For an affine root $\al+n$, define $U_{\al+n}$ to be the subgroup consisting of the elements $u_\al(xt^n)$ with $x\in\calO$. Set
\begin{equation*}
H_{\al,n}=\{x\in\calA_T\mid \al(x)=n\},
\qquad
H_{\al,>m}=\{x\in\calA_T\mid \al(x)>m\}.
\end{equation*}
We use analogous notation for the closed half-space $H_{\al,\geq m}$ and the strips $H_{\al,m<,<n}$, $H_{\al,m\leq,<n}$, and $H_{\al,m\leq,\leq n}$.
Write $\calH$ for the set of all affine root hyperplanes
$H_{\alpha,n}$, with $\alpha\in\Phi$ and $n\in\mathbb Z$.

Let $\Omega_0$ be the unique alcove contained in the $B$-positive Weyl chamber and having the origin as a vertex. This precise choice of base alcove will be used in Lemma~2.2. Let $W_{\mathrm{aff}}$ be the affine Weyl group generated by the reflections in the walls of $\Omega_0$. We have $W\cong W_{\mathrm{aff}}\rtimes\Omega$, where $\Omega$ is the subgroup of $W$ stabilizing $\Omega_0$. The Bruhat order on $W_{\mathrm{aff}}$ extends to $W$ by
\begin{equation*}
w_1\tau_1\geq w_2\tau_2
\quad\Longleftrightarrow\quad
w_1\geq w_2\ \text{ and }\ \tau_1=\tau_2,
\end{equation*}
for $w_i\in W_{\mathrm{aff}}$ and $\tau_i\in\Omega$.

\subsection{Affine flag variety notation}
Let $\calG$ be a group scheme over $k[[t]]$. Its loop functor $L\calG$ and positive loop functor $L^+\calG$ send a $k$-algebra $R$ to $\calG(R((t)))$ and $\calG(R[[t]])$, respectively. We denote by
\begin{equation*}
\operatorname{Fl}_{\calG}:=L\calG/L^+\calG
\end{equation*}
the \'etale sheafification of the presheaf quotient.

For every facet ${\bf f}\subset\overline{\Omega_0}$, Bruhat--Tits theory associates a parahoric group scheme $P_{\bf f}$ over $\calO$ whose $\calO$-points form the parahoric subgroup associated with ${\bf f}$ \cite{BT72}. The corresponding affine flag variety is represented by an ind-projective scheme over $k$.

Let ${\bf f}$ and ${\bf f}'$ be two facets of $\overline{\Omega_0}$. There is a natural action of $L^+P_{{\bf f}'}$ on $LG/L^+P_{\bf f}$. Let $W_{\bf f}$ and $W_{{\bf f}'}$ be the subgroups of $W$ generated by the reflections in the walls containing ${\bf f}$ and ${\bf f}'$, respectively, and let ${}_{{\bf f}'}W^{\bf f}$ denote the set of minimal-length representatives of $W_{{\bf f}'}\backslash W/W_{\bf f}$. The orbit of $w\in{}_{{\bf f}'}W^{\bf f}$ and its closure satisfy
\begin{equation} \label{closure relation}
\overline{L^+P_{{\bf f}'}wL^+P_{\bf f}/L^+P_{\bf f}}
=\bigsqcup_{\substack{v\in{}_{{\bf f}'}W^{\bf f}\\v\leq w}}
L^+P_{{\bf f}'}vL^+P_{\bf f}/L^+P_{\bf f}.
\end{equation}
See \cite[Proposition 2.8]{Ri13}. We denote by $Gr_{G,\mu}$ the orbit $L^+Gt^\mu L^+G/L^+G$ and by $Gr_{G,\leq\mu}$ its closure. These are the open and closed Schubert varieties in the affine Grassmannian $Gr_G$. We write $I$ for the positive loop group of the Iwahori group scheme of $G$ and $K=L^+G$.

\subsection{Retractions in $\calB(G,F)$}
We refer to \cite[Section 5.1]{Hai24} for the material in this subsection.

Fix a standard parabolic subgroup $P=MN$ of $G$, where $M$ is the standard Levi subgroup and $N$ is the unipotent radical. Let $K_P=LNL^+M$ and $I_P=LNI_M$, where $I_M=I\cap LM$. Choose an $M$-central, sufficiently $N$-dominant cocharacter $\nu\in X_*(T)$; thus
\begin{equation*}
\langle\al,\nu\rangle=0\quad(\al\in\Phi_M),
\qquad
\langle\be,\nu\rangle\gg0\quad(\be\in\Phi^+\setminus\Phi_M^+).
\end{equation*}
Write $t_\nu$ for geometric translation by $+\nu$. Thus
$t_\nu=t^{-\nu}$ in the Tits convention fixed above. Set
\[
\De_{I_P}=t_\nu\Omega_0=t^{-\nu}\Omega_0.
\]
This center lies far in the positive $N$-directions. Its image in the apartment for $M$ is the
base alcove. When $P=B$, the center $\De_U$ is far in the positive
Weyl chamber. This is the convention of \cite[Section 2.2 and
Lemma 5.3]{Hai24}. The resulting retraction is independent of the
sufficiently $N$-dominant choice in the sense made precise below.

For every alcove $\De\subset\calA_T$, there is an alcove-centered retraction
\begin{equation*}
\rho_{\De,\calA}:\calB(G,F)\longrightarrow\calA_T.
\end{equation*}
It restricts to the identity on $\calA_T$ and preserves distance from $\De$.

\begin{prop}[\cite{Hai24}, Lemma 5.3]
Let $\calY\subset\calB(G,F)$ be bounded. Then there is a bounded subgroup $I_{\calY}\subset I_P$ such that:
\begin{enumerate}
\item $\calY\subset\bigcup_{g\in I_{\calY}}g^{-1}\calA$;
\item the retractions $\rho_{t_\nu\Omega_0,\calA}$ agree on $\calY$ for every $M$-central and sufficiently $N$-dominant $\nu$, where the required degree of dominance depends on $\calY$.
\end{enumerate}
Consequently, their common value defines a retraction
\begin{equation*}
\rho_{I_P,\calA}:\calB(G,F)\longrightarrow\calA.
\end{equation*}
Moreover, $\rho_{I_P,\calA}^{-1}(\De)=I_P\De$ for every alcove $\De\subset\calA$.
\end{prop}

When $P=B$ and $U$ is the unipotent radical of $B$, the parabolic retraction $\rho_{I_P,\calA}$ is the retraction $\rho_{U,\calA}$ of \cite{BT72}. Before Section~7 we use only $\rho_{U,\calA}$; in Section~7 we use $\rho_{I_P,\calA}$ to prove the Levi branching rule.

\begin{lem}
Let $\Delta$ be an alcove of $\calB(G,F)$ having the origin as a vertex. Then
\begin{equation*}
\rho_{U,\calA}(\Delta)=\rho_{\Omega_0,\calA}(\Delta).
\end{equation*}
\end{lem}
\begin{proof}
One can argue by induction on the distance between $\Delta$ and
$\Omega_0$, so we may assume that there is an alcove
$\Delta'\subset\calA$ adjacent to $\Delta$, with a common
codimension $1$ face $F$, such that
$d(\Omega_0,\Delta')<d(\Omega_0,\Delta)$.
Let $\Delta''$ be the unique alcove in $\calA$ adjacent to
$\Delta'$ and having $F$ as a facet.
By \cite[Section 6.1]{HKM}, we have
\[
\rho_{U,\calA}(\Delta)
=\rho_{\Omega_0,\calA}(\Delta)
=\Delta'',
\]
and the lemma follows.
\end{proof}

\section{Galleries and alcove walks}

In this section we review galleries and alcove walks, define the raising and lowering operators in our conventions, compare them with those of \cite{GL05}, and prove the combinatorial lemmas used in Sections~4--6.

\subsection{Galleries and alcove walks}
\begin{dfn}
  A gallery in $\calB(G,F)$ is a sequence of facets
\begin{equation*}
\Gamma_0' \subset \Gamma_0 \supset \Gamma_1' \subset \Gamma_1
\supset \dots \subset \Gamma_s \supset \Gamma_{s+1}',
\end{equation*}
in $\calB(G,F)$, with the following properties:
\begin{enumerate}
    \item $\Gamma_0'$ and $\Gamma_{s+1}'$ are vertices.
    \item Each $\Gamma_i'$ is a codimension $1$ facet of both $\Gamma_{i-1}$ and $\Gamma_i$ for $1 \leq i \leq s$.
    \item All $\Gamma_i$'s are of the same dimension.
\end{enumerate}  

When all the $\Gamma_i$ are alcoves, this is called a vertex-to-vertex alcove walk.
\end{dfn}
\begin{rmk}
If all the $\Gamma_i$ are alcoves and we remove the initial and terminal vertices, the resulting chain is called an alcove-to-alcove walk. The alcove walks in \cite{Hai24} are alcove-to-alcove walks beginning at the base alcove $\Omega_0$. Unless otherwise stated, the term alcove walk in this paper means a vertex-to-vertex alcove walk.
\end{rmk}

For the rest of the paper, fix $\mu\in X_*^+(T)$ and put $\mu^*:=-w_0\mu$, where $w_0$ is the longest element of $W_0$. Define
\begin{equation*}
H_{\mu^*}:=
\bigcap_{\substack{\al\in\Phi\\\langle\al,\mu^*\rangle=0}}H_{\al,0}.
\end{equation*}
Thus $H_{\mu^*}$ is the intersection of the root hyperplanes whose reflections fix $\mu^*$.

\begin{dfn}[\cite{GL05}, Definition 11 and Lemma 4]
    A gallery
    \begin{equation*}
    \ga_{\mu}:\Ga_0'\subset \Ga_0\supset \Ga_1'
    \subset\cdots\subset \Ga_s\supset \Ga_{s+1}'
    \end{equation*}
    is a minimal gallery from $0$ to $\mu^*$ if it is contained in $H_{\mu^*}$ and:
    \begin{enumerate}
        \item $\Ga_0'=0$ and $\Ga_{s+1}'=\mu^*$;
        \item $\dim(\Ga_i)=\dim H_{\mu^*}$ for $0\leq i\leq s$;
        \item if $\calH_i$ is the set of hyperplanes containing $\Ga_i'$ but not $\Ga_i$, then
        \begin{equation*}
        \bigsqcup_{0\leq i\leq s}\calH_i=\calM(0,\mu^*),
        \end{equation*}
        where $\calM(0,\mu^*)$ is the set of affine root hyperplanes separating $0$ and $\mu^*$, with $0$ in the closed half-space and $\mu^*$ in the complementary open half-space.
    \end{enumerate}
\end{dfn}

Fix a minimal gallery from $0$ to $\mu^*$,
\[
\ga_\mu:\Upsilon_0'\subset\Upsilon_0\supset\Upsilon_1'
\subset\cdots\subset\Upsilon_s\supset\Upsilon_{s+1}',
\qquad
\Upsilon_0'=0,\quad\Upsilon_{s+1}'=\mu^*.
\]
We construct its alcove refinement as in
\cite[Section~4, proof of Lemma~4]{GL05}. For a face $F$ and
an alcove $C$, let $\operatorname{proj}_F(C)$ be the unique
alcove containing $F$ in its closure at minimal gallery
distance from $C$.

Set $i_0=0$. Since $\mu^*$ is dominant,
$\Upsilon_0\subset\overline{\Omega_0}$.
If $s\geq1$, choose a minimal alcove walk
$(\Omega_0,\Omega_1,\ldots,\Omega_{i_1})$ from $\Omega_0$
to $\Omega_{i_1}=\operatorname{proj}_{\Upsilon_1}(\Omega_0)$.
Inductively, for $2\leq k\leq s$, suppose that the construction has reached
$\Omega_{i_{k-1}}$, with
$\Upsilon_{k-1}\subset\overline{\Omega_{i_{k-1}}}$.
Choose a minimal alcove walk
$(\Omega_{i_{k-1}},\Omega_{i_{k-1}+1},\ldots,\Omega_{i_k})$
such that
\[
\Omega_{i_k}
=\operatorname{proj}_{\Upsilon_k}(\Omega_{i_{k-1}}),
\qquad
i_k-i_{k-1}=d(\Omega_{i_{k-1}},\Omega_{i_k}).
\]
For each $1\leq k\leq s$, both endpoints of the $k$th block
contain $\Upsilon_k'$ in their closures.
The residue of $\Upsilon_k'$ is convex, so
$\Upsilon_k'\subset\overline{\Omega_i}$ for
$i_{k-1}\leq i\leq i_k$. Put $p=i_s$.

By the wall calculation in the proof of
\cite[Lemma~4]{GL05}, the $k$th block crosses precisely the
walls containing $\Upsilon_k'$ but not $\Upsilon_k$.
Minimality of $\ga_\mu$ makes these sets pairwise disjoint.
Consequently, the concatenation is a minimal alcove walk,
and $\Omega_p=\operatorname{proj}_{\{\mu^*\}}(\Omega_0)$.
With $\Omega_i'$ the common panel of $\Omega_{i-1}$ and
$\Omega_i$, we obtain
\[
\de_\mu:\{0\}\subset\Omega_0\supset\Omega_1'
\subset\Omega_1\supset\cdots\supset\Omega_p'
\subset\Omega_p\supset\{\mu^*\}.
\]
In particular, the construction gives
\[
\Upsilon_k\subset\overline{\Omega_{i_k}}
\quad(0\leq k\leq s),\qquad
\Upsilon_k'\subset\overline{\Omega_i}
\quad(i_{k-1}\leq i\leq i_k,\ 1\leq k\leq s).
\]

We fix the choice of $\ga_\mu$ and $\de_\mu$ for the rest of this article.

Let $m_k$ and $n_k$ be the types of $\Upsilon_k'$ and
$\Upsilon_k$, respectively. For a face type $t$, let
$F_t\subset\overline{\Omega_0}$ be the face of type $t$, let
$W_t$ be its pointwise stabilizer in $W_{\mathrm{aff}}$, and
let $w_t$ be the longest element of $W_t$.
The refinement above is obtained by choosing reduced expressions
\begin{equation}\label{eq:type-blocks}
s_{t_{i_{k-1}+1}}\cdots s_{t_{i_k}}
=w_{m_k}w_{n_k}\qquad(1\leq k\leq s),
\end{equation}
where $t_i$ is the type of $\Omega_i'$.
Here $w_{m_k}w_{n_k}$ is the minimal-length representative of
the maximal coset in $W_{m_k}/W_{n_k}$; see
\cite[Section~4, Proposition~2]{GL05} and the discussion following it.
Concatenating these expressions gives
$\Omega_i=s_{t_1}\cdots s_{t_i}\Omega_0$ and fixes a reduced
expression
\begin{equation}\label{eq:extended-word}
(t_{\mu^*})_0=s_{t_1}\cdots s_{t_p}\,\omega_\mu,
\qquad \omega_\mu\in\Omega.
\end{equation}
As usual, $y_0$ denotes the minimal-length element of $yW_0$.
The factor $\omega_\mu$ records the terminal vertex type and
the component of the affine Grassmannian.

\begin{dfn}
\label{def:fixed-gallery-types}
We fix $\ga_\mu$, its refinement $\de_\mu$, and the indices
$i_k$ throughout. Write $t_{\ga_\mu}$ and $t_{\de_\mu}$ for
their chains of face types. Let $\Ga(\ga_\mu)$ be the set of
galleries in $\calA_T$ of type $t_{\ga_\mu}$ starting at $0$,
and let $\De(\de_\mu)$ be the set of alcove walks in $\calA_T$
of type $t_{\de_\mu}$ starting at $0$. Write
$\Ga(\lambda,\ga_\mu)$ and $\De(\lambda,\de_\mu)$ for the
subsets whose endpoint is $-\lambda$.
\end{dfn}

\begin{dfn}[Skeleton]
\label{def:skeleton}
For $\de\in\De(\de_\mu)$, write
$\de:\De_0'\subset\De_0\supset\De_1'\subset\cdots
\subset\De_p\supset\De_{p+1}'$.
For $0\leq k\leq s$, let $\Ga_k$ be the unique face of
$\overline{\De_{i_k}}$ of type $n_k$. For $1\leq k\leq s$,
let $\Ga_k'$ be the unique face of
$\overline{\De_{i_{k-1}}}$ of type $m_k$, and put
$\Ga_0'=\De_0'$ and $\Ga_{s+1}'=\De_{p+1}'$.
The skeleton of $\de$ is
\[
\operatorname{sk}(\de):
\Ga_0'\subset\Ga_0\supset\Ga_1'\subset\cdots
\subset\Ga_s\supset\Ga_{s+1}'.
\]
\end{dfn}

The construction in Definition~\ref{def:skeleton} gives a
well-defined map
\[
\operatorname{sk}:\De(\lambda,\de_\mu)
\longrightarrow\Ga(\lambda,\ga_\mu).
\]
For $1\leq k\leq s$, the small face $\Ga_k'$ is the unique
face of type $m_k$ contained in every
$\overline{\De_i}$ with $i_{k-1}\leq i\leq i_k$.
In particular, $\operatorname{sk}(\de_\mu)=\ga_\mu$.

\begin{rmk}
One may also consider alcove walks associated with arbitrary reduced
expressions of $(t_{\mu^*})_0$. However, intersecting the closed alcoves
of an arbitrary minimal walk with $H_{\mu^*}$ need not give a minimal
gallery; the intersections need not even form a gallery.
For example, in type $B_2$, take coordinates $(x,y)$ with
$\al_1(x,y)=x$ and $\al_2(x,y)=y-x$, and put $\mu^*=(0,4)$.
Then $H_{\mu^*}=\{x=0\}$ and
\[
\overline{\Omega_0}
=\operatorname{conv}\{(0,0),(1/2,1/2),(0,1)\},\qquad
(t_{\mu^*})_0=s_0s_1s_0s_2s_1s_0s_2s_1s_0,
\]
where $s_0(x,y)=(1-y,1-x)$, $s_1(x,y)=(-x,y)$ and
$s_2(x,y)=(y,x)$. Write $\Omega_i$ for the alcove obtained by
applying the first $i$ letters of this word to $\Omega_0$.
The walk crosses each of the nine walls $x+y=j$, $y=j$ and
$y-x=j$, $1\leq j\leq3$, exactly once. These are precisely the
walls separating $\Omega_0$ from $\mu^*$, so the walk is minimal.
Nevertheless,
\[
\overline{\Omega_6}
=\operatorname{conv}\{(1,2),(1/2,5/2),(1,3)\},\qquad
\overline{\Omega_6}\cap H_{\mu^*}=\varnothing.
\]
The one-dimensional intersections, after removing repetitions, are
$\{0\}\times[0,1]$, $\{0\}\times[1,2]$ and $\{0\}\times[3,4]$.
The missing edge $\{0\}\times[2,3]$ prevents them from forming a gallery.

\begin{figure}[htbp]
\centering
\begin{tikzpicture}[scale=1.35]
  \begin{scope}
    \clip (-.35,-.12) rectangle (1.45,4.15);
    \foreach \j in {-1,...,5} {
      \draw[thin,lightgray!80] (-.35,\j) -- (1.45,\j);
      \draw[thin,lightgray!80] (-.35,\j-.35) -- (1.45,\j+1.45);
      \draw[thin,lightgray!80] (-.35,\j+.35) -- (1.45,\j-1.45);
    }
    \foreach \j in {0,1} {
      \draw[thin,lightgray!80] (\j,-.12) -- (\j,4.15);
    }
    \fill[blue!10] (1,2) -- (.5,2.5) -- (1,3) -- cycle;
    \foreach \a/\b/\c in {
      {0,0}/{.5,.5}/{0,1},
      {1,1}/{.5,.5}/{0,1},
      {1,1}/{.5,1.5}/{0,1},
      {0,2}/{.5,1.5}/{0,1},
      {0,2}/{.5,1.5}/{1,2},
      {0,2}/{.5,2.5}/{1,2},
      {1,3}/{.5,2.5}/{1,2},
      {1,3}/{.5,2.5}/{0,3},
      {1,3}/{.5,3.5}/{0,3},
      {0,4}/{.5,3.5}/{0,3}} {
      \draw[thick] (\a) -- (\b) -- (\c) -- cycle;
    }
  \end{scope}
  \draw[red,ultra thick] (0,0) -- (0,2);
  \draw[red,ultra thick] (0,3) -- (0,4);
  \draw[red,thick,dashed] (0,2) -- (0,3);
  \foreach \i/\x/\y in {
    0/.1667/.5,1/.5/.8333,2/.5/1.1667,3/.1667/1.5,
    4/.5/1.8333,5/.5/2.1667,6/.8333/2.5,
    7/.5/2.8333,8/.5/3.1667,9/.1667/3.5} {
    \coordinate (walk-\i) at (\x,\y);
  }
  \foreach \i [evaluate=\i as \j using int(\i+1)] in {0,...,8} {
    \draw[blue!75!black,thick,-{Stealth[length=1.3mm]},
      shorten <=2pt,shorten >=2pt] (walk-\i) -- (walk-\j);
  }
  \foreach \i in {0,...,9} {
    \node[font=\scriptsize,fill=white,inner sep=.5pt] at (walk-\i) {$\i$};
  }
  \foreach \j in {0,...,4} {
    \fill[red] (0,\j) circle (1.2pt);
    \node[left,font=\small] at (-.04,\j) {$\j$};
  }
  \node[above,font=\small] at (0,4.17) {$H_{\mu^*}=\{x=0\}$};
  \node[above right,font=\small] at (0,4) {$\mu^*$};
  \node[left,font=\small,text=red!75!black] at (-.12,2.5) {missing edge};
  \draw[blue!75!black,thick,-{Stealth[length=1.5mm]}]
    (1.8,2.9) -- (2.2,2.9);
  \node[right,font=\small] at (2.25,2.9) {minimal alcove walk};
  \draw[red,ultra thick] (1.8,2.5) -- (2.2,2.5);
  \node[right,font=\small] at (2.25,2.5) {one-dimensional intersections};
  \node[anchor=west,font=\small] at (1.8,2.1) {The number $i$ labels $\Omega_i$.};
\end{tikzpicture}
\caption{A minimal alcove walk in type $B_2$ whose intersections with $H_{\mu^*}$ do not form a gallery.}
\label{fig:B2-intersection-counterexample}
\end{figure}
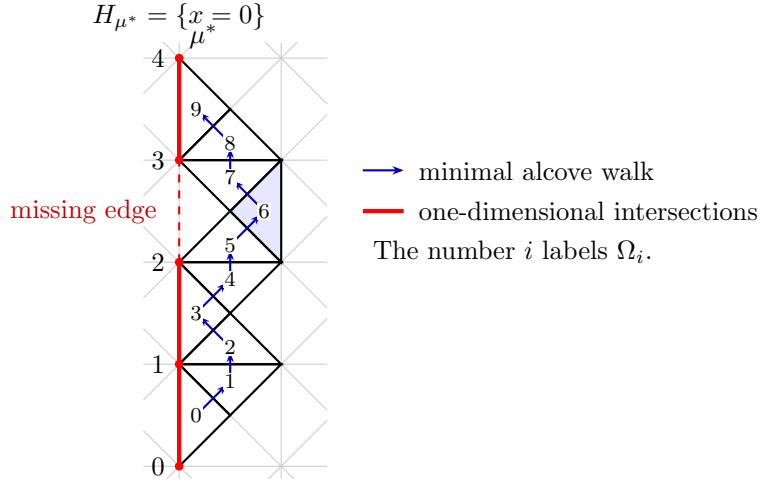

For maximal-dimensional positively folded walks, restricting to the
chosen type of $\de_\mu$ entails no loss of information.
Indeed, Lemma~\ref{lem:braid-existence} changes any reduced expression
of $(t_{\mu^*})_0$ to the one fixed by $\de_\mu$, preserving the initial
and terminal alcoves and the dimension. On maximal-dimensional walks this construction is independent of the
sequence of braid moves, and reversing the moves gives its inverse.
\end{rmk}

By \cite[Proposition 2]{GL05}, there is a bijection
\begin{equation*}
W_{m_0}\times^{W_{n_0}}W_{m_1}\times^{W_{n_1}}\cdots
\times^{W_{n_{s-1}}}W_{m_s}/W_{n_s}
\xrightarrow{\ \sim\ }\Ga(\ga_\mu).
\end{equation*}
It sends $(v_0,\ldots,v_s)$ to the gallery whose $i$th large face is
\begin{equation*}
\Ga_i=v_0v_1\cdots v_i F_{n_i}.
\end{equation*}

We next define $\De$-positive foldings. Let $\De$ be an alcove of $\calA_T$, and let
\begin{equation*}
\ga:\Ga_0'\subset\Ga_0\supset\cdots\supset\Ga_{s+1}'
\end{equation*}
be represented by $(v_0,\ldots,v_s)$, where each $v_i$ is the minimal-length representative of its class in $W_{m_i}/W_{n_i}$. Let $\tau_i$ be the minimal-length representative of the unique maximal element of $W_{m_i}/W_{n_i}$. Then $\ga_\mu$ is represented by $(1,\tau_1,\ldots,\tau_s)$.

For later comparison with alcove walks, we record the crossing convention first.
\begin{dfn}[\cite{Hai24}, Definitions~4.1 and~4.2]
\label{def:alcove-crossings-foldings}
Let
\begin{equation*}
\de:\De_0'\subset\De_0\supset\cdots\supset\De_{p+1}'
\end{equation*}
be an alcove walk, and let $H$ be the wall spanned by $\De_k'$. A crossing at $\De_k'$ is $\De$-positive (respectively, $\De$-negative) if $\De_k$ and $\De$ lie on opposite (respectively, the same) sides of $H$. A folding at $\De_k'$ is $\De$-positive (respectively, $\De$-negative) if $\De_{k-1}=\De_k$ and this common alcove and $\De$ lie on opposite (respectively, the same) sides of $H$.
\end{dfn}

\begin{prop}[\cite{GL05}, Lemma 5]\label{prop:fold-reflections}
For $1\leq i\leq s$, if $v_i\neq\tau_i$, the gallery $\ga$ is folded at the smaller face $\Ga_i'$. Let $\ga^{i_1}$ be represented by
\begin{equation*}
(v_0,\ldots,v_{i-1},\tau_i,\ldots,\tau_s)
\end{equation*}
and let $\ga^i$ be represented by
\begin{equation*}
(v_0,\ldots,v_i,\tau_{i+1},\ldots,\tau_s).
\end{equation*}
Write the large faces of $\ga^{i_1}$ as $\Ga_j$ and those of $\ga^i$ as $\widetilde{\Ga_j}$ for $j\geq i$. Then there are positive roots $\be_1,\ldots,\be_m$ and integers $a_1,\ldots,a_m$ such that
\begin{equation}\label{eq:fold-reflections}
\widetilde{\Ga_j}
=s_{\be_m,a_m}\cdots s_{\be_1,a_1}\Ga_j
\qquad (j\geq i).
\end{equation}
The reflections may be taken in walls containing $\Ga_i'$.
Writing $g_i=v_0\cdots v_{i-1}$, we have
\[
s_{\be_m,a_m}\cdots s_{\be_1,a_1}
\in g_i v_iW_{n_i}\tau_i^{-1}g_i^{-1}.
\]
Thus the product fixes $\Ga_i'$ and sends
$g_i\tau_iF_{n_i}$ to $g_iv_iF_{n_i}$.
\end{prop}

\begin{dfn}\label{def:positive-folding}
With the notation of Proposition~\ref{prop:fold-reflections},
the folding at $\Ga_i'$ is $\De$-positive if there exist
positive roots $\be_1,\ldots,\be_m$ and integers
$a_1,\ldots,a_m$ such that
\[
E^{(0)}=g_i\tau_iF_{n_i},\qquad
E^{(j)}=s_{\be_j,a_j}E^{(j-1)}
\quad (1\leq j\leq m),\qquad
E^{(m)}=g_iv_iF_{n_i},
\]
and, for every $1\leq j\leq m$, the wall $H_{\be_j,a_j}$
contains $\Ga_i'$ and separates $E^{(j)}$ from $\De$.
A gallery is $\De$-positively folded if every folding is
$\De$-positive. Write $\Ga_\De^+(\ga_\mu)$ for these galleries.
\end{dfn}

For alcove-to-alcove walks, this convention agrees with \cite[Definitions 4.1 and 4.2]{Hai24}: a fold is positive precisely when its wall separates the folded alcove from the center $\De$. For $\De=\De_U$, our retraction is taken in the direction opposite to that used in \cite[Definition 16]{GL05}. Thus, after applying the involution $x\mapsto-x$ to the apartment, our positive foldings are exactly the positive foldings of Gaussent--Littelmann. 

We use only the centers $\De_U$ and $\De_{I_P}$. Write $\Ga^+(\ga_\mu)$ and $\Ga_P^+(\ga_\mu)$ for the corresponding positively folded galleries. Unless a center is specified, positive means $\De_U$-positive. We use analogous notation for alcove walks.

We next define the dimension of a gallery. For $\De=\De_U$ and $\De=\De_{I_P}$, Sections~4 and~7 identify this combinatorial statistic with the dimension of the corresponding Bia{\l}ynicki--Birula cell in the Bott--Samelson variety.

\begin{dfn}
Let $\ga\in\Ga(\ga_\mu)$. A hyperplane $H\in\calH$ is a $\De$-load-bearing wall for $\ga$ at the smaller face $\Ga_i'$ if $\Ga_i'\subset H$ and $H$ separates $\Ga_i$ from $\De$. Without a subscript, load-bearing means $\De_U$-load-bearing.
\end{dfn}
The dimension of a gallery is defined by the following count.

\begin{dfn}
Let $\ga\in\Ga(\ga_\mu)$. Define
\begin{equation*}
\dim_{\De}(\ga):=
\#\bigl\{(H,\Ga_i')\bigm|
0\leq i\leq s,\ 
H\text{ is }\De\text{-load-bearing for }\ga\text{ at }\Ga_i'\bigr\}.
\end{equation*}
We abbreviate $\dim_{\De_U}(\ga)$ to $\dim(\ga)$.
\end{dfn}

\begin{rmk}\label{rmk:walk-dimension}
Let $\de$ be a positively folded walk of type $t_{\de_\mu}$ with
first alcove $w\Omega_0$. Prepend a minimal walk from $\Omega_0$
to $w\Omega_0$. The resulting alcove-to-alcove walk $\de'$ has
type $w(t_{\mu^*})_0$ by Lemma~\ref{lem:minimal-translation}.
Its dimension in \cite[Definition 4.7]{Hai24} equals $\dim(\de)$:
the load-bearing walls at the initial vertex of $\de$ are exactly
the walls crossed in the prepended segment. Section~4 uses this
equality to compare our vertex-to-vertex walks with Haines' model.
\end{rmk}

\subsection{Root operators and comparison with \cite{GL05}}
\label{subsec:gallery-root-operators}

We define the Gaussent--Littelmann root operators for galleries in
$\calA_T$ of arbitrary type, including vertex-to-vertex alcove walks.
Our galleries start at $0$ and end at $-\lambda$, and translations
act by the Tits convention. Our retraction direction is opposite
to that in \cite{GL05}; these choices explain the differences from
\cite[Definitions~19 and~20]{GL05}.

\begin{dfn}\label{def:gallery-root-operators}
Let
\begin{equation*}
\ga:\Ga_0'\subset\Ga_0\supset\cdots\supset\Ga_{r+1}'
\end{equation*}
be a gallery in $\calA_T$ with $\Ga_0'=0$ and
$\Ga_{r+1}'=-\lambda$, where $\lambda\in X_*(T)$.
No fixed-type or positive-folding assumption is imposed.
For a simple root $\al$, put
\[
m=m_\al(\ga)
:=\max\{n\in\mathbb Z\mid\Ga_i'\subset H_{\al,n}
\text{ for some }i\}.
\]
The two endpoints give $m\geq0$ and
$m\geq-\langle\al,\lambda\rangle$. The cut indices below exist
because a gallery meets each intermediate affine wall in a small face.

\begin{enumerate}
\item If $m=0$, then $e_\al^{GL}(\ga)$ is undefined.
If $m>0$, let
\[
k=\min\{i\mid\Ga_i'\subset H_{\al,m}\},\qquad
j=\max\{i<k\mid\Ga_i'\subset H_{\al,m-1}\}.
\]
On the large faces, define
\begin{equation}\label{eq:raw-raising}
\bigl(e_\al^{GL}(\ga)\bigr)_i=
\begin{cases}
\Ga_i,&i<j,\\
s_{\al,m-1}\Ga_i,&j\leq i<k,\\
\Ga_i-\al^\vee,&i\geq k.
\end{cases}
\end{equation}
The smaller faces are transformed by the same maps. At the two
joining faces, the maps agree on $H_{\al,m-1}$ and $H_{\al,m}$,
respectively. The tail translation is $t^{\al^\vee}$, so the
endpoint is $-(\lambda+\al^\vee)$.

\item If $m=-\langle\al,\lambda\rangle$, then
$f_\al^{GL}(\ga)$ is undefined. Otherwise, define independently
of the preceding indices
\[
k=\max\{i\mid\Ga_i'\subset H_{\al,m}\},\qquad
j=\min\{i>k\mid\Ga_i'\subset H_{\al,m-1}\}.
\]
On the large faces, define
\begin{equation}\label{eq:raw-lowering}
\bigl(f_\al^{GL}(\ga)\bigr)_i=
\begin{cases}
\Ga_i,&i<k,\\
s_{\al,m}\Ga_i,&k\leq i<j,\\
\Ga_i+\al^\vee,&i\geq j.
\end{cases}
\end{equation}
Again the smaller faces are transformed by the same maps, which
agree at the joining faces. The tail translation is
$t^{-\al^\vee}$, so the endpoint is $-(\lambda-\al^\vee)$.
\end{enumerate}
When defined, these operations preserve the gallery type.
\end{dfn}

\begin{rmk}
Let $\iota:\calA_T\to\calA_T$ be the involution $x\mapsto-x$. Under $\iota$, the maximal hyperplane level $m$ becomes the minimal level $-m$ used in \cite[Definitions 19 and 20]{GL05}. Conjugating our foldings and tail translations by $\iota$ gives exactly the operators of \cite{GL05}.
\end{rmk}

We now specialize to the fixed type $\ga_\mu$ and recall the consequences from \cite{GL05}. For every $\ga\in\Ga^+(\lambda,\ga_\mu)$,
\begin{equation*}
\dim(\ga)\leq\langle\rho,\lambda+\mu\rangle.
\end{equation*}
Equality defines an LS gallery. Write $\Ga_{LS}^+(\ga_\mu)$ for all LS galleries and $\Ga_{LS}^+(\lambda,\ga_\mu)$ for those ending at $-\lambda$. For $\ga\in\Ga_{LS}^+(\lambda,\ga_\mu)$, whenever the indicated operator is defined,
\begin{align*}
e_\al^{GL}(\ga)&\in\Ga_{LS}^+(\lambda+\al^\vee,\ga_\mu),\\
f_\al^{GL}(\ga)&\in\Ga_{LS}^+(\lambda-\al^\vee,\ga_\mu).
\end{align*}
This explicit compatibility is the input used in Sections~5 and~6 to transport the Gaussent--Littelmann crystal structure to the maximal-dimensional positively folded alcove walks.
    
\subsection{Combinatorial lemmas}

We now record the combinatorial facts used to construct the two Bott--Samelson varieties in Section~4 and to compare their strata in Sections~5 and~6.

Let $W_{\mu^*}$ be the stabilizer of $\mu^*$ in $W_0$, and let $W^{\mu^*}$ be the set of minimal-length representatives of $W_0/W_{\mu^*}$. Every $w\in W_0$ has a unique length-additive decomposition
\begin{equation*}
w=w^{\mu^*}w_{\mu^*},
\qquad
w^{\mu^*}\in W^{\mu^*},\quad w_{\mu^*}\in W_{\mu^*}.
\end{equation*}
Let $w_{0,\mu^*}$ be the longest element of $W_{\mu^*}$.
We use $W_\mu$ and $w_{0,\mu}$ for the corresponding stabilizer
and its longest element.

\begin{lem}\label{lem:minimal-translation}
Let $(t_{\mu^*})_0$ be the unique minimal-length element of the coset $t_{\mu^*}W_0$. Then
\begin{equation}\label{eq:minimal-translation}
(t_{\mu^*})_0=t_{\mu^*}w_0w_{0,\mu}.
\end{equation}
Consequently, for every $v\in W_0$,
\begin{equation*}
\ell\bigl(v(t_{\mu^*})_0\bigr)
=\ell(v)+\ell\bigl((t_{\mu^*})_0\bigr).
\end{equation*}
For $w\in W_{\mu^*}$,
\begin{equation}\label{eq:move-stabilizer}
w(t_{\mu^*})_0=(t_{\mu^*})_0\mu(w),
\qquad
\mu(w):=w_{0,\mu}w_0ww_0w_{0,\mu}.
\end{equation}
\end{lem}
\begin{proof}
The minimality assertion and \eqref{eq:minimal-translation} are the standard description of the minimal element in this double coset; compare \cite[Section 5.3]{Hai24}. Minimality on the left gives the length formula. Since $w(\mu^*)=\mu^*$, conjugating $w$ by \eqref{eq:minimal-translation} gives
\begin{align*}
(t_{\mu^*})_0^{-1}w(t_{\mu^*})_0
&=w_{0,\mu}w_0t_{-\mu^*}wt_{\mu^*}w_0w_{0,\mu}\\
&=w_{0,\mu}w_0ww_0w_{0,\mu}
=\mu(w),
\end{align*}
which proves \eqref{eq:move-stabilizer}.
\end{proof}
\begin{lem}\label{lem:partial-conjugation}
For $1\leq k\leq s$, put
\begin{equation*}
x_k=s_{t_1}\cdots s_{t_{i_k}}.
\end{equation*}
If $w\in W_{\mu^*}$, then the concatenation of a reduced word for $w$ with the displayed word for $x_k$ is reduced, and
\begin{equation}\label{eq:partial-reduced}
\ell(wx_k)=\ell(w)+i_k.
\end{equation}
Moreover, there is a unique element $\mu_k(w)\in W_{n_k}$ such that
\begin{equation}\label{eq:partial-conjugation}
wx_k=x_k\mu_k(w),
\qquad
\ell(\mu_k(w))=\ell(w).
\end{equation}
At the final index,
\begin{equation*}
\mu_s(w)=\omega_\mu\mu(w)\omega_\mu^{-1}.
\end{equation*}
\end{lem}

\begin{proof}
The alcoves of $\de_\mu$ lie in the dominant Weyl chamber. Hence, for every simple reflection $s_\al\in W_{\mu^*}$ and every initial subword $s_{t_1}\cdots s_{t_j}$,
\begin{equation*}
\ell\bigl(s_\al s_{t_1}\cdots s_{t_j}\bigr)
=1+\ell\bigl(s_{t_1}\cdots s_{t_j}\bigr).
\end{equation*}
Thus $x_k$ is a minimal representative for the left cosets of
$W_{\mu^*}$. The parabolic length-additivity criterion gives
\eqref{eq:partial-reduced} for every $w\in W_{\mu^*}$.

The face $\Upsilon_k\subset\overline{\Omega_{i_k}}$ has full
dimension in $H_{\mu^*}$ and is fixed by $W_{\mu^*}$.
Since $\Omega_{i_k}$ lies in the dominant chamber, its chamber
in the transverse spherical residue is the dominant chamber
for $W_{\mu^*}$. Thus $H_{\alpha,0}$ is a wall of
$\Omega_{i_k}$ for every simple root $\alpha$ fixing $\mu^*$.
Transporting these walls to $\Omega_0$ shows that conjugation
by $x_k$ sends the simple generators of $W_{\mu^*}$ to the
simple generators of $W_{n_k}$. Therefore
\begin{equation*}
\mu_k(w):=x_k^{-1}wx_k\in W_{n_k},
\end{equation*}
and this conjugation preserves the length of reduced words in these generators. This proves \eqref{eq:partial-conjugation}. For $k=s$, use $x_s=(t_{\mu^*})_0\omega_\mu^{-1}$
in \eqref{eq:move-stabilizer}. This gives the displayed terminal
formula. Conjugation by $\omega_\mu$ permutes affine simple
reflections and preserves length.
\end{proof}
\begin{lem}\label{lem:coset-alternative}
Let $w\in W^{\mu^*}$, $v\in W_{\mu^*}$, and let $s_\al$ be a simple reflection such that
\begin{equation*}
\ell(s_\al wv)>\ell(wv).
\end{equation*}
Exactly one of the following alternatives occurs:
\begin{enumerate}
\item $(s_\al wv)^{\mu^*}=w$, in which case $w(\mu^*)\in H_{\al,0}$;
\item $(s_\al wv)^{\mu^*}=s_\al w$ and $(s_\al wv)_{\mu^*}=v$.
\end{enumerate}
\end{lem}
\begin{proof}
First, $\ell(s_\al w)=\ell(w)+1$; otherwise,
\begin{equation*}
\ell(s_\al wv)\leq\ell(s_\al w)+\ell(v)
<\ell(w)+\ell(v)=\ell(wv),
\end{equation*}
a contradiction. If $s_\al w\in W^{\mu^*}$, length additivity gives the second alternative.

Suppose $s_\al w\notin W^{\mu^*}$. By the standard criterion for minimal coset representatives (compare \cite[Lemma 3.1]{HH17}), there is a simple reflection $s_\be\in W_{\mu^*}$ such that
\begin{equation*}
s_\al ws_\be=w.
\end{equation*}
Thus $s_\al=ws_\be w^{-1}$ fixes $w(\mu^*)$, and the first alternative follows.
\end{proof}

\section{Bott--Samelson varieties and the alcove-walk model}
\subsection{The two Bott--Samelson varieties}

Using the types fixed in Section~3, define
\begin{align*}
X_1
&:=L^+G\times^I L^+P_{t_1}\times^I\cdots
   \times^I L^+P_{t_p}/I,\\
X_2
&:=L^+G\times^{L^+P_{n_0}}L^+P_{m_1}
   \times^{L^+P_{n_1}}\cdots
   \times^{L^+P_{n_{s-1}}}L^+P_{m_s}/L^+P_{n_s}.
\end{align*}
Here $P_{t_i}$, $P_{m_j}$, and $P_{n_j}$ are the parahoric group schemes attached to facets of the indicated types. Recall that $i_s=p$. Multiplication within the blocks singled out by the indices $i_j$ defines a morphism
\begin{equation}\label{eq:skeleton-morphism}
\begin{aligned}
\pi:X_1&\longrightarrow X_2,\\
(g_0,g_1,\ldots,g_p)
&\longmapsto
\bigl(g_0,\,
g_1\cdots g_{i_1},\,
g_{i_1+1}\cdots g_{i_2},\ldots,
g_{i_{s-1}+1}\cdots g_{i_s}\bigr).
\end{aligned}
\end{equation}
The contracted-product relations show that $\pi$ is well defined. A $k$-point of $X_1$ (respectively, $X_2$) represents a gallery in $\calB(G,F)$ of type $t_{\de_\mu}$ (respectively, $t_{\ga_\mu}$), and, by Definition~\ref{def:skeleton}, $\pi$ sends a gallery to its skeleton.

\begin{rmk}
In the split case, the Bruhat--Tits interpretation of the $k$-points of these Bott--Samelson varieties does not require $k$ to be perfect; see \cite[Section 2]{Hai25}.
\end{rmk}

Multiplication also defines endpoint morphisms
\begin{equation*}
\begin{aligned}
f:X_1&\longrightarrow Gr_{G,\leq\mu},
& (g_0,\ldots,g_p)&\longmapsto g_0\cdots g_p\dot\omega_\mu K,\\
g:X_2&\longrightarrow Gr_{G,\leq\mu},
& (h_0,\ldots,h_s)&\longmapsto h_0\cdots h_s\dot\omega_\mu K.
\end{aligned}
\end{equation*}
Here $\dot\omega_\mu$ is a fixed representative in $N_G(T)(F)$.
The terminal parahoric is contained in
$\dot\omega_\mu K\dot\omega_\mu^{-1}$, so these formulas are well
defined. On $k$-points, each morphism sends a gallery to its
marked endpoint. In particular, $f=g\circ\pi$.

\begin{prop}
The morphism $g$ is a resolution of singularities of $Gr_{G,\leq\mu}$, and its restriction over the open Schubert cell is an isomorphism:
\begin{equation*}
g^{-1}(Gr_{G,\mu})\xrightarrow{\ \sim\ }Gr_{G,\mu}.
\end{equation*}
\end{prop}
\begin{proof}
The contracted product $X_2$ is an iterated fibration by projective homogeneous spaces and is therefore smooth and projective. The multiplication morphism associated with the reduced gallery type is the standard Demazure resolution of the Schubert variety $Gr_{G,\leq\mu}$. It is birational and is an isomorphism over the open Schubert cell. This is the affine-flag resolution constructed by Richarz; see \cite[Section 2]{Ri13}, together with \cite[Proposition 3.7.4]{dHL} for the description over the open cells.
\end{proof}

\begin{rmk}
The formulation above is the one used in this paper; compare \cite[Remark 3.27]{CvdHS25} and \cite[Remark 14]{GL05}.
\end{rmk}

We next analyze $f$. The Iwahori-orbit decomposition of its first factor gives
\begin{align*}
X_1
&=\bigsqcup_{w\in W_0}
IwI\times^I L^+P_{t_1}\times^I\cdots\times^I L^+P_{t_p}/I\\
&=:\bigsqcup_{w\in W_0}C(w).
\end{align*}
By Lemma~\ref{lem:minimal-translation}, if a reduced expression for $w$ is followed by
\begin{equation*}
s_{t_1}s_{t_2}\cdots s_{t_p}\omega_\mu=(t_{\mu^*})_0,
\end{equation*}
the resulting extended expression for $w(t_{\mu^*})_0$ is reduced. Multiplication, followed by right multiplication by $\dot\omega_\mu$, gives a restriction of a Demazure morphism
\begin{equation*}
d_w:C(w)\longrightarrow\overline{Iw(t_{\mu^*})_0I/I}.
\end{equation*}
It is birational and is an isomorphism over
$Iw(t_{\mu^*})_0I/I$; see \cite[Proposition 8.8]{PR08}.

\subsection{Retraction strata on the Bott--Samelson varieties}

Assume first that $k$ is algebraically closed. Applying $\rho_{U,\calA}$ to every face of a gallery gives maps
\begin{equation*}
\rho^1_{U,\calA}:X_1(k)\longrightarrow\De(\de_\mu),
\qquad
\rho^2_{U,\calA}:X_2(k)\longrightarrow\Ga(\ga_\mu).
\end{equation*}
The fibers are the retraction strata, denoted $C_{X_1}(\de)$ and $C_{X_2}(\ga)$.

Over an arbitrary field $k$, use a generic dominant cocharacter to define the $\bbG_m$-action on the first factor of $X_i$, with limits as $c\to0$. Indeed, conjugation by this cocharacter contracts $LU$; the opposite choice contracts $LU^-$ and describes the opposite retraction. Its attractor functor is
\begin{equation*}
\begin{aligned}
X_i^+:\operatorname{AffSch}_k^{\mathrm{opp}}
&\longrightarrow\operatorname{Sets},\\
R&\longmapsto
\operatorname{Hom}^{\bbG_m}_R(\bbA_R^{1,+},X_{i,R}).
\end{aligned}
\end{equation*}
The attractor is represented by a scheme, with evaluation morphisms
\[
\operatorname{ev}_0:X_i^+\longrightarrow X_i^{\bbG_m},
\qquad
\operatorname{ev}_1:X_i^+\longrightarrow X_i;
\]
see \cite[Section~1]{Ri19}. The fixed points are the combinatorial walks $\de\in\De(\de_\mu)$
for $X_1$ and the galleries $\ga\in\Ga(\ga_\mu)$ for $X_2$.
Since this limit agrees with facewise $U$-retraction, evaluation
at $1$ identifies the fibers over these fixed points with the
retraction strata:
\begin{equation}\label{eq:attractor-retraction-cells}
\operatorname{ev}_0^{-1}(\de)
\xrightarrow[\operatorname{ev}_1]{\sim}C_{X_1}(\de),
\qquad
\operatorname{ev}_0^{-1}(\ga)
\xrightarrow[\operatorname{ev}_1]{\sim}C_{X_2}(\ga).
\end{equation}
These identifications give the strata their locally closed scheme
structures in $X_i$. Thus $X_i^+$ is the disjoint union of the
attracting cells, whose images under $\operatorname{ev}_1$ give a partition
of $X_i$ into locally closed strata.
The root-group coordinates supplied by the fixed split pinning are defined over $k$. Therefore the Bia{\l}ynicki--Birula paving and each of its retraction strata are defined over $k$. For $X_2$, the strata are indexed by $\Ga(\ga_\mu)$, and \cite[Lemma 13]{GL05} gives
\begin{equation*}
C_{X_2}(\ga)\cong\bbA_k^{\dim(\ga)}
\end{equation*}
for every positively folded gallery $\ga$. The analogous statement holds for $X_1$ and its alcove-walk strata. See also \cite{BB73,Ri19}.

The morphism $\pi$ is compatible with the retractions. More precisely, for $\de\in\De(\de_\mu)$,
\begin{equation}\label{eq:pi-retraction-strata}
\pi\bigl(C_{X_1}(\de)\bigr)
\subset C_{X_2}\bigl(\operatorname{sk}(\de)\bigr).
\end{equation}

Let $q:LG/I\to LG/K$ be the natural projection. For $w\in W_0$, write $w=w^{\mu^*}w_{\mu^*}$. The double-coset calculation in \cite[Section 5.6]{HKP10} gives
\begin{equation}\label{eq:minimal-coset-representative}
(wt_{\mu^*})_0=(w^{\mu^*}t_{\mu^*})_0
=w^{\mu^*}(t_{\mu^*})_0.
\end{equation}
For $w\in W^{\mu^*}$, the restriction
\begin{equation*}
q:Iw(t_{\mu^*})_0I/I
\xrightarrow{\ \sim\ }
Iw(t_{\mu^*})_0K/K
\end{equation*}
is an isomorphism by \cite[Proposition 4.4]{Hai25}. Moreover, \cite[Lemma 7.2 and Proposition 10.1]{Hai25} identify
\begin{equation}\label{eq:semi-infinite-projection}
\bigsqcup_{u\in W_0}
\bigl(LUt^\lambda uI/I\cap Iw(t_{\mu^*})_0I/I\bigr)
\xrightarrow{\ \sim\ }
LUt^\lambda K/K\cap Iw(t_{\mu^*})_0K/K.
\end{equation}

For general $w\in W_0$, equation \eqref{eq:minimal-coset-representative} identifies the image of $Iw(t_{\mu^*})_0I/I$ with the $I$-orbit $Iw^{\mu^*}(t_{\mu^*})_0K/K$ in $LG/K$. Composing $q$ with the inverse of the isomorphism for $w^{\mu^*}$ defines
\begin{equation*}
f_w:Iw(t_{\mu^*})_0I/I
\longrightarrow Iw^{\mu^*}(t_{\mu^*})_0I/I.
\end{equation*}
Write $w=w^{\mu^*}v$. Lemma~\ref{lem:minimal-translation} gives
the length-additive factorization
\[
w(t_{\mu^*})_0=w^{\mu^*}(t_{\mu^*})_0\mu(v),
\qquad \mu(v)\in W_0,\quad \ell(\mu(v))=\ell(v).
\]
Choose root-group coordinates for this reduced expression.
The coordinates belonging to the final factor $\mu(v)$ form
an affine space of dimension $\ell(v)$. Right multiplication
by this factor does not change the coset modulo $K$.
Consequently these coordinates identify $f_w$ with the
projection
\[
Iw^{\mu^*}(t_{\mu^*})_0I/I\times\mathbb A^{\ell(v)}
\longrightarrow Iw^{\mu^*}(t_{\mu^*})_0I/I.
\]
The next theorem describes precisely how the semi-infinite intersections pull back. It is used below to compare the fibers over different initial alcoves and to define the sets $B_\mu^v(\lambda)$.

\begin{thm}
For every $w\in W_0$,
\begin{equation}\label{eq:fw-preimage}
f_w^{-1}\!\left(
\bigsqcup_{u\in W_0}
LUt^\lambda uI/I\cap Iw^{\mu^*}(t_{\mu^*})_0I/I
\right)
=
\bigsqcup_{u\in W_0}
LUt^\lambda uI/I\cap Iw(t_{\mu^*})_0I/I.
\end{equation}
\end{thm}
\begin{proof}
Let $q:LG/I\to LG/K$ be the projection. The finite flag
decomposition of $K/I$ gives
\[
q^{-1}(LUt^\lambda K/K)
=\bigsqcup_{u\in W_0}LUt^\lambda uI/I.
\]
By definition, $f_w$ is the restriction of $q$ followed by
the inverse of the isomorphism
\[
Iw^{\mu^*}(t_{\mu^*})_0I/I
\xrightarrow{\sim}Iw^{\mu^*}(t_{\mu^*})_0K/K.
\]
Both sides of \eqref{eq:fw-preimage} are therefore
\[
Iw(t_{\mu^*})_0I/I\ \cap\ q^{-1}(LUt^\lambda K/K).
\]
All schemes carry their reduced structures, as stipulated in
Section~2, so this also gives the asserted equality of
locally closed subvarieties.
\end{proof}

The projection $q:LG/I\to LG/K$ restricts to a morphism
\begin{equation}\label{eq:pi-prime}
\pi'(w):
\overline{Iw(t_{\mu^*})_0I/I}
\longrightarrow
\overline{Iw^{\mu^*}(t_{\mu^*})_0K/K}.
\end{equation}
Indeed, by the closure relation \eqref{closure relation}, each stratum $IxI/I$ in the source maps to $Ix^0K/K$, where $x^0$ is the minimal element of $xW_0$. If $x\leq w(t_{\mu^*})_0$, then
\begin{equation*}
x^0\leq\bigl(w(t_{\mu^*})_0\bigr)^0
=w^{\mu^*}(t_{\mu^*})_0;
\end{equation*}
compare \cite[Lemma 3.1]{HH17}. This proves \eqref{eq:pi-prime}.

Compose the Demazure morphism $d_w$ with \eqref{eq:pi-prime} and set
\begin{equation*}
D_w:=\pi'(w)\circ d_w:
C(w)\longrightarrow
\overline{Iw^{\mu^*}(t_{\mu^*})_0K/K}.
\end{equation*}

\begin{prop}
There is an equality
\begin{equation*}
D_w^{-1}\bigl(Iw^{\mu^*}(t_{\mu^*})_0K/K\bigr)
=d_w^{-1}\bigl(Iw(t_{\mu^*})_0I/I\bigr).
\end{equation*}
We denote this common open subvariety by $C^0(w)$.
\end{prop}
\begin{proof}
A point of $C(w)$ is represented by an alcove-to-alcove gallery of type $w(t_{\mu^*})_0$. Under the retraction $\rho_{\Omega_0,\calA}$, its first $\ell(w)+1$ alcoves form a minimal walk from $\Omega_0$ to $w\Omega_0$; the remaining part has type $\de_\mu$.

Both sides of the asserted equality may be tested after applying the retraction to the vertices of the gallery. A point lies over the open $I$-cell precisely when the retracted gallery has no folding after $w\Omega_0$. If a folding occurs, the Weyl distance from $w\Omega_0$ to the terminal alcove is strictly less than $p=\ell((t_{\mu^*})_0)$. Hence its terminal vertex cannot be $w\mu^*$, which is the retracted vertex characterizing the indicated open $I$-orbit in $LG/K$. Conversely, an unfolded minimal gallery has terminal vertex $w\mu^*$ and lies over both open cells. This proves the equality.
\end{proof}
We now explain why Theorem~4.4 and Proposition~4.5 are needed. The variety $X_1$ parametrizes galleries of type $\de_\mu$, whereas $X_2$ parametrizes their skeletons of type $\ga_\mu$. The map $\pi$ forgets the intermediate alcoves, and $g$ identifies $g^{-1}(Gr_{G,\mu})$ with $Gr_{G,\mu}$. These facts allow us to compare the retraction stratifications over the MV intersection
\begin{equation*}
Z_{\lambda,\mu}:=LUt^\lambda K/K\cap Kt^\mu K/K.
\end{equation*}

The Iwahori decomposition of the open Schubert cell gives
\begin{equation}\label{eq:MV-Iwahori-decomposition}
Z_{\lambda,\mu}
=\bigsqcup_{w\in W^{\mu^*}}
\bigl(LUt^\lambda K/K
\cap Iw(t_{\mu^*})_0K/K\bigr).
\end{equation}
The index $w$ has an intrinsic meaning: it is determined by the Iwahori orbit containing the generic point of the corresponding irreducible component.

Fix $w\in W^{\mu^*}$ and choose $w'\in W_0$ with $(w')^{\mu^*}=w$. Theorem~4.4 and Proposition~4.5 identify the inverse image under $D_{w'}$ with
\begin{equation*}
\bigl(LUt^\lambda K/K\cap Iw(t_{\mu^*})_0K/K\bigr)
\times\bbA^{\ell(w'_{\mu^*})}.
\end{equation*}
Via $C^0(w')\cong Iw'(t_{\mu^*})_0I/I$, this inverse image decomposes as
\begin{equation*}
\bigsqcup_{u\in W_0}
\bigl(LUt^\lambda uI/I\cap Iw'(t_{\mu^*})_0I/I\bigr).
\end{equation*}
By \cite[Lemma 5.4]{Hai24}, its irreducible components of dimension
\begin{equation*}
\langle\rho,\lambda+\mu\rangle+\ell(w'_{\mu^*})
\end{equation*}
are indexed by the positively folded alcove-to-alcove walks of type $w'(t_{\mu^*})_0$ that begin at $\Omega_0$ and end at $t_{-\lambda}u\Omega_0$ for some $u\in W_0$ and have dimension $\langle\rho,\lambda+\mu\rangle+\ell(w'_{\mu^*})$. Denote this set by $M_\mu(w',\lambda)$. Thus there are natural bijections
\begin{align*}
&\operatorname{Irr}^{\langle\rho,\lambda+\mu\rangle}
\bigl(LUt^\lambda K/K\cap Iw(t_{\mu^*})_0K/K\bigr)\\
&\quad\longleftrightarrow
\operatorname{Irr}^{\langle\rho,\lambda+\mu\rangle
  +\ell(w'_{\mu^*})}
\bigl(D_{w'}^{-1}(LUt^\lambda K/K
  \cap Iw(t_{\mu^*})_0K/K)\bigr)\\
&\quad\longleftrightarrow M_\mu(w',\lambda).
\end{align*}
In other words, we have the following proposition:

\begin{prop}\label{prop:MV-walk-bijection}
For every $v\in W_{\mu^*}$, there is a natural bijection
\begin{equation}\label{eq:MV-walk-bijection}
\operatorname{Irr}^{\langle\rho,\lambda+\mu\rangle}
(Z_{\lambda,\mu})
\longleftrightarrow
\bigsqcup_{w\in W^{\mu^*}}M_\mu(wv,\lambda).
\end{equation}
\end{prop}

The right-hand side is therefore an alcove-walk model for the MV cycles in $Z_{\lambda,\mu}$.
\begin{rmk}\label{rmk:haines-alcove-model}
When $v=e$, the model on the right of \eqref{eq:MV-walk-bijection} recovers Haines' model \cite[Theorem 5.15]{Hai24} in the Borel case.
\end{rmk}

\begin{lem}\label{lem:initial-minimal-segment}
Let $x\in W_0$. If $\delta\in M_\mu(x,\lambda)$ has type $x(t_{\mu^*})_0$, then its first $\ell(x)+1$ alcoves form a minimal alcove walk from $\Omega_0$ to $x\Omega_0$.
\end{lem}
\begin{proof}
This follows from positive folding and Lemma~2.2: no wall in the initial segment may fold toward the base alcove, so the segment is minimal.
\end{proof}

For $w\in W^{\mu^*}$, $v\in W_{\mu^*}$, and
$\de\in M_\mu(wv,\lambda)$, the \emph{initial minimal segment}
means the first $\ell(wv)$ steps, from $\Omega_0$ to $wv\Omega_0$,
as in Lemma~\ref{lem:initial-minimal-segment}.
Remove these steps, retaining $wv\Omega_0$ as the first alcove,
and mark $0$ and $-\lambda$ as the initial and terminal vertices.
The remainder is a vertex-to-vertex alcove walk of type $\de_\mu$;
Remark~\ref{rmk:walk-dimension} identifies its dimension with that
of $\de$. This gives the following definition.

\begin{dfn}\label{def:alcove-model}
For $w\in W^{\mu^*}$ and $v\in W_{\mu^*}$, let $B_\mu(wv,\lambda)$ be the set of positively folded vertex-to-vertex alcove walks $\de$ of type $t_{\de_\mu}$ whose first alcove is $wv\Omega_0$, whose endpoint is $-\lambda$, and whose dimension is
\begin{equation*}
\dim(\de)=\langle\rho,\lambda+\mu\rangle+\ell(v).
\end{equation*}
Define
\begin{equation*}
B_\mu^v(\lambda)
:=\bigsqcup_{w\in W^{\mu^*}}B_\mu(wv,\lambda),
\qquad
B_\mu^v:=\bigsqcup_\lambda B_\mu^v(\lambda),
\end{equation*}
where the second union is over the cocharacters $\lambda$ in the convex hull of $W_0\mu$ such that $\mu-\lambda\in\mathbb Z\Phi^\vee$.
\end{dfn}

\begin{prop}\label{prop:MV-vertex-walk-bijection}
For every $v\in W_{\mu^*}$, there is a natural bijection
\[
\operatorname{Irr}^{\langle\rho,\lambda+\mu\rangle}
(Z_{\lambda,\mu})
\longleftrightarrow B_\mu^v(\lambda).
\]
\end{prop}
\begin{proof}
Removing the initial minimal segment from $\Omega_0$ to $wv\Omega_0$ gives a bijection
\[
\bigsqcup_{w\in W^{\mu^*}}M_\mu(wv,\lambda)
\longleftrightarrow B_\mu^v(\lambda),
\]
whose inverse prepends that segment. Compose this with the bijection in Proposition~\ref{prop:MV-walk-bijection}.
\end{proof}

\section{Bijections between different models}\label{sec:bijections}

The geometric comparisons in Section~4 identify several combinatorial
parametrizations of the MV cycles in
\[
Z_{\lambda,\mu}:=LUt^\lambda K/K\cap Kt^\mu K/K.
\]
We first show that the skeleton map identifies the maximal-dimensional
alcove walks in $B_\mu^v(\lambda)$ with LS galleries. This bijection will
be used in Section~6 to transport the Gaussent--Littelmann crystal
structure to $B_\mu^v$. We then compare the different choices of
$v\in W_{\mu^*}$ and give a braid-move construction of
these comparison maps.

\subsection{The skeleton bijection}

\begin{prop}\label{prop:skeleton-bijection}
For every $v\in W_{\mu^*}$ and every weight $\lambda$, the map
\[
\begin{aligned}
\operatorname{PR}_v:B_\mu^v(\lambda)
&\longrightarrow \Gamma_{LS}^+(\lambda,\ga_\mu),\\
\de&\longmapsto\operatorname{sk}(\de)
\end{aligned}
\]
is a bijection.
\end{prop}

\begin{proof}
Fix $w\in W^{\mu^*}$ and $\de\in B_\mu(wv,\lambda)$, and put
$\ga=\operatorname{sk}(\de)$. As in Lemma~\ref{lem:initial-minimal-segment} and Definition~\ref{def:alcove-model},
prepend the minimal initial segment from $\Omega_0$ to $wv\Omega_0$.
This gives the alcove-to-alcove walk of type $wv(t_{\mu^*})_0$ used
in the parametrization of Section~4. We keep the notation $\de$
for its stratum in $X_1$.

Use the identifications of Proposition~4.5 and write
\[
\begin{aligned}
Y_{wv}&:=Iwv(t_{\mu^*})_0I/I, &
Y_w&:=Iw(t_{\mu^*})_0K/K,\\
S_\de&:=C_{X_1}(\de)\cap Y_{wv},&
T_\ga&:=C_{X_2}(\ga)\cap Y_w.
\end{aligned}
\]
The endpoint maps identify $Y_{wv}$ with $C^0(wv)$ and $Y_w$
with its inverse image in the open Schubert cell of $X_2$.
The projection
\[
\pi_{wv}:Y_{wv}\longrightarrow Y_w
\]
has fibers isomorphic to $\mathbb A^{\ell(v)}$. The compatibility
of $\pi$ with retractions, established in Section~4.2, gives
\[
\pi_{wv}(S_\de)\subseteq T_\ga.
\]
The source is a nonempty open subset of the retraction stratum
indexed by $\de$. Hence
\[
\dim S_\de=\dim\de
=\langle\rho,\lambda+\mu\rangle+\ell(v).
\]
The dimension-of-fibers inequality now gives
\begin{align*}
\dim\ga
&\geq\dim T_\ga\\
&\geq\dim\pi_{wv}(S_\de)\\
&\geq\dim S_\de-\ell(v)\\
&=\langle\rho,\lambda+\mu\rangle.
\end{align*}
The stratum $T_\ga$ meets the open locus of minimal galleries,
so $\ga$ is positively folded by \cite[Theorem 4]{GL05}.
The reverse inequality is therefore the gallery dimension bound
recalled in Section~3.2. Thus $\ga$ is an LS gallery.
Moreover, equality holds throughout, so $\pi_{wv}(S_\de)$ is dense
in $T_\ga$.

We now compare the component parametrizations. Equation~\eqref{eq:MV-walk-bijection} identifies the walks in $B_\mu^v(\lambda)$
with the irreducible components of $Z_{\lambda,\mu}$. The LS-gallery
parametrization identifies $\ga$ with the closure of its retraction
stratum in the same intersection. The density just proved shows that
these two parametrizations associate the same component to $\de$
and $\operatorname{sk}(\de)$. Consequently, the skeleton map is their
composition and is a bijection.
\end{proof}

\begin{dfn}\label{def:initial-direction}
For an LS gallery
\[
\ga:\Ga_0'\subset\Ga_0\supset\cdots\supset\Ga_{s+1}',
\]
let $t(\ga)$ be the unique representative in $W^{\mu^*}$ such that
\[
\Ga_0=t(\ga)\Upsilon_0.
\]
Here the stabilizer of $\Upsilon_0$ is $W_{\mu^*}$, so its translates
are indexed by $W_0/W_{\mu^*}$. For $w\in W^{\mu^*}$, set
\[
\Gamma_{LS,w}^+(\lambda,\ga_\mu)
:=\{\ga\in\Gamma_{LS}^+(\lambda,\ga_\mu)\mid t(\ga)=w\}.
\]
The restriction of Proposition~\ref{prop:skeleton-bijection} is
\[
pr_{wv}:B_\mu(wv,\lambda)
\xrightarrow{\ \sim\ }\Gamma_{LS,w}^+(\lambda,\ga_\mu).
\]
For $\ga$ in this target, write
\[
\ga^{wv}:=pr_{wv}^{-1}(\ga).
\]
\end{dfn}

Thus $pr_{wv}$ records the initial direction, whereas
$\operatorname{PR}_v$ is the union of these maps over $w$.
The maximal-dimension hypothesis in
Proposition~\ref{prop:skeleton-bijection} is essential, as the
second example below shows.

\begin{ex}\label{ex:A2-skeleton}
    \begin{figure}[htbp]
    \centering
    \begin{minipage}[t]{0.46\textwidth}
    \centering
    % ======= 第一张图 (左) =======
    \begin{tikzpicture}[scale=1.2] % 稍微缩小了 scale (从 1.5 改为 1.2) 以防止超出页面宽度
        % 定义正三角形高度 h = sqrt(3)/2
        \def\h{0.866025}
        
        % 保持原来的裁剪空间不变
        \clip (-0.7, -0.5) rectangle (4.2, 4);

        % 绘制完美对齐的 A2 仿射公寓网格 (细灰线)
        \begin{scope}[gray, thin]
            \foreach \i in {-6,...,6} {
                % 水平线
                \draw (-6, \i*\h) -- (6, \i*\h);
                % 60度方向斜线
                \draw (\i - 3, -6*\h) -- (\i + 3, 6*\h);
                % 120度方向斜线
                \draw (\i + 3, -6*\h) -- (\i - 3, 6*\h);
            }
        \end{scope}

        % 绘制 5 个嵌入的三角形 (Alcoves)
        \begin{scope}[ultra thick, line join=round]
            
            % 1. 内部共享边 (黑色)
            \draw (1, 2*\h) -- (1.5, 3*\h);   
            \draw (1, 2*\h) -- (2, 2*\h);     
            \draw (1.5, 1*\h) -- (2, 2*\h);   
            \draw (1.5, 1*\h) -- (2.5, 1*\h); 

            % 2. 顶部和右侧的外部边界 (黑色)
            \draw (0.5, 3*\h) -- (1.5, 3*\h) -- (2.5, 1*\h) -- (2, 0);

            % 3. 偏左的三条外部边界 (红色)
            \draw[red] (0.5, 3*\h) -- (1, 2*\h) -- (1.5, 1*\h) -- (2, 0);

            % --- 放置标签符号 ---
            \node[above left, font=\large] at (0.5, 3*\h) {$\mu^*$};
            \node[above, font=\large] at (2, 3.5*\h) {$\calC_+$};
            \node[left, font=\large, red] at (1.25, 1.5*\h) {$\ga_{\mu}$};
            \node[above left, font=\large] at (0,4*\h) {$H_{\al_1}$};
            \node[above right, font=\large] at (3.5, 3*\h) {$H_{\al_2}$};
            \node[font=\large] at (1.75, 0.5*\h) {1};
            \node[font=\large] at (2.25, 0.5*\h) {2};
            \node[font=\large] at (2,1*\h) {0};

        \end{scope}

    \end{tikzpicture}
    \caption{The minimal alcove walk $\de_\mu$ and its skeleton $\ga_\mu$.}
    \label{fig:A2-minimal-skeleton}
    \end{minipage}%
    \hspace{0.5cm} % <-- 在这里调整左右两张图的间距
    \begin{minipage}[t]{0.46\textwidth}
    \centering
    % ======= 第二张图 (右) =======
    \begin{tikzpicture}[scale=1.2]
        % 定义正三角形高度 h = sqrt(3)/2
        \def\h{0.866025}
        
        % 保持原来的裁剪空间不变
        \clip (-0.7, -0.5) rectangle (4.2, 4);

        % 绘制完美对齐的 A2 仿射公寓网格 (细灰线)
        \begin{scope}[gray, thin]
            \foreach \i in {-6,...,6} {
                % 水平线
                \draw (-6, \i*\h) -- (6, \i*\h);
                % 60度方向斜线
                \draw (\i - 3, -6*\h) -- (\i + 3, 6*\h);
                % 120度方向斜线
                \draw (\i + 3, -6*\h) -- (\i - 3, 6*\h);
            }
        \end{scope}

        % 绘制 5 个嵌入的三角形 (Alcoves)
        \begin{scope}[ultra thick, line join=round]
            
            % 1. 内部共享边 (黑色)  
            \draw (1, 2*\h) -- (2, 2*\h);     
            \draw (1.5, 1*\h) -- (2, 2*\h);   
            \draw (1.5, 1*\h) -- (2.5, 1*\h); 
            \draw (0.5, 1*\h) -- (1.5, 1*\h);
            \draw (1.15, 1.925*\h) -- (1.85, 1.925*\h);

            % 2. 顶部和右侧的外部边界 (黑色)
            \draw (2, 2*\h) -- (2.5, 1*\h) -- (2, 0);

            % 3. 偏左的三条外部边界 (红色)
            \draw[red] (1, 2*\h) -- (1.5, 1*\h) -- (2, 0);

            \draw[red] (0.5, 1*\h) -- (1, 2*\h);

            % --- 放置标签符号 ---
            \node[above left, font=\large] at (0.5, 3*\h) {$\mu^*$};
            \node[above, font=\large] at (2, 3.5*\h) {$\calC_+$};
            \node[left, font=\large, red] at (1.25, 1.5*\h) {$\ga$};
            \node[above left, font=\large] at (0,4*\h) {$H_{\al_1}$};
            \node[above right, font=\large] at (3.5, 3*\h) {$H_{\al_2}$};
            \node[font=\large] at (1.75, 0.5*\h) {1};
            \node[font=\large] at (2.25, 0.5*\h) {2};
            \node[font=\large] at (2,1*\h) {0};

        \end{scope}

    \end{tikzpicture}
    \caption{The LS gallery $\ga$ and its maximal-dimensional alcove-walk lift $\de$.}
    \label{fig:A2-LS-lift}
    \end{minipage}
\end{figure}

Let $G=\mathrm{SL}_3$, let $T$ be the diagonal torus, and let $B$
be the Borel subgroup of upper triangular matrices. Take the
positive simple roots and dominant coweights
\[
\begin{aligned}
\al_1&=(1,-1,0),& \al_2&=(0,1,-1),\\
\mu&=3w_1=(2,-1,-1),&
\mu^*&=3w_2=(1,1,-2),
\end{aligned}
\]
where $w_1,w_2$ are the fundamental coweights. The panels of the
base alcove are labeled by their affine simple reflections.
The left diagram shows the chosen minimal walk of type
$(t_{\mu^*})_0=s_0s_1s_2s_0$ and its skeleton $\ga_\mu$ in red.

The right diagram shows the LS gallery
\[
\ga=e_{\al_1}^{GL}e_{\al_2}^{GL}(\ga_\mu)
\]
and its maximal-dimensional lift $\de\in B_\mu^e$.
Starting from the base alcove, the walk has a negative crossing
of type $s_0$, a negative crossing of type $s_1$, a positive
fold of type $s_2$, and a positive crossing of type $s_0$.
Counting load-bearing walls gives
\[
\dim\de=\dim\ga=2.
\]
\end{ex}

\begin{rmk}\label{rmk:B2-nonmaximal}
The skeleton map need not be injective on all positively folded
walks. In particular, it does not in general identify the entire
set $\calP_\mu^{\De_U}(\lambda)$ of
\cite[Definition 5.9]{Hai24} with $\Ga^+(\lambda,\ga_\mu)$.
These sets index two stratifications of the same MV intersection,
but only their maximal-dimensional parts are identified by
Proposition~\ref{prop:skeleton-bijection}.

For an example, let $G=\mathrm{SO}_5$ be the adjoint group of type
$B_2$. In the usual coordinate realizations of the character and
cocharacter lattices, take
\[
\begin{aligned}
\al_1&=\epsilon_1,&
\al_2&=-\epsilon_1+\epsilon_2,\\
\al_1^\vee&=2\epsilon_1,&
\al_2^\vee&=-\epsilon_1+\epsilon_2.
\end{aligned}
\]
Here $\epsilon_1,\epsilon_2$ are the coordinate vectors, and
\[
\mu=\mu^*=2w_2=2\epsilon_2.
\]
Choose $\de_\mu$ of type $(t_{\mu^*})_0=s_0s_1s_0$.
Its skeleton $\ga_\mu$ lies in $H_{\al_1,0}$. The LS gallery
\[
\ga=e_{\al_2}^{GL}e_{\al_1}^{GL}e_{\al_2}^{GL}(\ga_\mu)
:\Ga_0'\subset\Ga_0\supset\Ga_1'\subset\Ga_1\supset\Ga_2'
\]
has $\Ga_0=\Ga_1$, the one-dimensional face of the base alcove
of type $s_1$, and $\dim\ga=3$.

Both walks in the diagrams below start at the base alcove.
The left walk $\de_1$ has a positive fold of type $s_0$, a
positive crossing of type $s_1$, and a positive fold of type $s_0$.
The right walk $\de_2$ has a negative crossing of type $s_0$,
a positive fold of type $s_1$, and a positive crossing of type $s_0$.
Thus
\[
\operatorname{sk}(\de_1)=\operatorname{sk}(\de_2)=\ga,
\qquad
\dim\de_1=3,\qquad \dim\de_2=2.
\]
Only $\de_1$ is the maximal-dimensional lift belonging to $B_\mu^e$.

\begin{figure}[htbp]
    \centering
    \begin{minipage}[t]{0.44\textwidth}
    \centering
    % ======= 第一张图 (左) =======
    \begin{tikzpicture}[scale=1.2]
        % 绘制所有细线条网格（作为背景和未加粗部分）
        \foreach \i in {0,...,4} {
            \draw[thin, lightgray!80] (\i,0) -- (\i,4);
            \draw[thin, lightgray!80] (0,\i) -- (4,\i);
        }
        
        % 绘制所有 X 形细对角线
        \foreach \x in {0,...,3} {
            \foreach \y in {0,...,3} {
                \draw[thin, lightgray!80] (\x,\y) -- (\x+1,\y+1);
                \draw[thin, lightgray!80] (\x,\y+1) -- (\x+1,\y);
            }
        }
        
        % 加粗菱形
        \draw[thick] (2,2) -- (1.5,2.5); 
        \draw[thick] (2,2) -- (2.5,2.5);
        \draw[thick] (1.5,2.5) -- (2,3); 
        \draw[thick] (2.5,2.5) -- (2,3); 
        \draw[thick] (2,2) -- (2,3);
        \draw[thick] (2.05,2.90) -- (2.45,2.5);
        \draw[thick] (1.55,2.5) -- (1.95,2.9);
        \draw[thick] (2.25,2.5) -- (1.75,2.5);
        \draw[thick] (1.75,2.5) -- (1.8,2.55);
        \draw[thick] (1.75,2.5) -- (1.8,2.45);

        \node[above ] at (2, 3.5) {$H_{\al_1}$};
        \node[above right] at (3.5,3.5) {$H_{\al_2}$};

    \end{tikzpicture}
    \caption{The maximal-dimensional lift $\de_1$ of $\ga$ in type $B_2$.}
    \label{fig:B2-maximal-lift}
    \end{minipage}%
    \hspace{1.5cm} % <-- 在这里调整两张图左右之间的距离
    \begin{minipage}[t]{0.44\textwidth}
    \centering
    % ======= 第二张图 (右) =======
    \begin{tikzpicture}[scale=1.2]
        % 绘制所有细线条网格（作为背景和未加粗部分）
        \foreach \i in {0,...,4} {
            \draw[thin, lightgray!80] (\i,0) -- (\i,4);
            \draw[thin, lightgray!80] (0,\i) -- (4,\i);
        }
        
        % 绘制所有 X 形细对角线
        \foreach \x in {0,...,3} {
            \foreach \y in {0,...,3} {
                \draw[thin, lightgray!80] (\x,\y) -- (\x+1,\y+1);
                \draw[thin, lightgray!80] (\x,\y+1) -- (\x+1,\y);
            }
        }
        
        % 加粗菱形
        
        \draw[thick] (2,2) -- (2.5,2.5);
        \draw[thick] (2,2) -- (2,3);
        \draw[thick] (2.5,2.5) -- (2,3); 
        \draw[thick] (2,3) -- (3,3);
        \draw[thick] (2.5,2.5) -- (3,3);
        \draw[thick] (2.1,2.95) -- (2.9,2.95);

        \node[above ] at (2, 3.5) {$H_{\al_1}$};
        \node[above right] at (3.5,3.5) {$H_{\al_2}$} ;
        
    \end{tikzpicture}
    \caption{The non-maximal-dimensional lift $\de_2$ of $\ga$ in type $B_2$.}
    \label{fig:B2-nonmaximal-lift}
    \end{minipage}
\end{figure}

\end{rmk}

\subsection{Changing the initial alcove}

The comparison of two choices of $v$ is already determined by the
skeleton bijections. We record this before discussing how to compute it.

\begin{prop}\label{prop:transition}
For $w\in W^{\mu^*}$ and $v,v'\in W_{\mu^*}$, define
\[
\begin{aligned}
T^w_{v\leftarrow v'}:B_\mu(wv',\lambda)
&\xrightarrow{\ \sim\ }B_\mu(wv,\lambda),\\
T^w_{v\leftarrow v'}&:=pr_{wv}^{-1}\circ pr_{wv'}.
\end{aligned}
\]
This is the unique bijection preserving the skeleton. It satisfies
\[
T^w_{v\leftarrow v}=\mathrm{id},
\qquad
T^w_{v\leftarrow v'}\circ T^w_{v'\leftarrow v''}
=T^w_{v\leftarrow v''}.
\]
For every $\de\in B_\mu(wv',\lambda)$,
\[
\dim T^w_{v\leftarrow v'}(\de)
=\dim\de+\ell(v)-\ell(v').
\]
\end{prop}
\begin{proof}
The assertions about bijections and composition follow from
Definition~\ref{def:initial-direction}. The dimension formula follows
from Definition~\ref{def:alcove-model}.
\end{proof}

We describe the construction of $T^w_{v\leftarrow e}$.
The general comparison follows by composition. We first record
the braid-move statement needed for independence of choices.

\begin{lem}\label{lem:braid-existence}
Fix the orientation determined by $\De_U$.
A positively folded walk of reduced type can be replaced by
a walk of any other reduced expression for the same element,
with the same initial and terminal alcoves and the same number
of folds and dimension.

Suppose, moreover, that a walk $\xi$ of reduced type $x$ from
$\Omega_0$ to $y\Omega_0$ satisfies
\[
\dim\xi=\dim(LUyI/I\cap IxI/I).
\]
Then any two sequences of braid moves connecting the same
initial and final reduced expressions give the same walk.
\end{lem}

\begin{proof}
The existence assertion is \cite[Proposition~4.33]{MST19}.
We prove the additional assertions over an algebraic closure
of $k$. Put
\[
Z(x,y)=LUyI/I\cap IxI/I,\qquad d=\dim Z(x,y).
\]
For a reduced expression $\mathbf i$ of $x$ and a point
$z\in IxI/I$, let $\widetilde\delta_{\mathbf i}(z)$ be
the unique minimal gallery of type $\mathbf i$ from
$\Omega_0$ to the alcove represented by $z$.
The endpoint morphism on the open Bott--Samelson variety
identifies it with $IxI/I$. Thus
\[
S_{\mathbf i}(\xi)
=\{z\in Z(x,y):
  \rho_{U,\calA}(\widetilde\delta_{\mathbf i}(z))=\xi\}
\]
is the walk stratum of \cite[Lemma~5.4]{Hai24}.
For each $\mathbf i$, these strata give a finite partition
of $Z(x,y)$ into irreducible locally closed subvarieties,
and $\dim S_{\mathbf i}(\xi)=\dim\xi$.

First we show that for a specific braid move $st\cdots\to ts\cdots$, there is a unique way to change the alcove walk correspondingly if it is of maximal dimension. Consider a braid block of length $m=m(s,t)$.
It is a reduced expression for the longest element of the
finite group $\langle s,t\rangle$. Fix its initial and
terminal retracted alcoves $A$ and $B$.
Let $R$ be the finite flag variety of their rank-two residue,
and put
\[
\begin{aligned}
\mathcal O_A&=\{D\in R:\operatorname{dist}(A,D)=m\},\\
\mathcal O_B&=\{D\in R:\rho_{U,\calA}(D)=B\},
\end{aligned}
\]
where $\operatorname{dist}$ is gallery distance.
The first is the open Bruhat cell of chambers opposite $A$.
The second is a Bruhat cell for the Borel fixing the
projection of the sector-center to this residue.
The endpoint variety for this block is
$ 
\mathcal O_A\cap\mathcal O_B.
$ 
It is nonempty whenever the block occurs, and is an open
subset of the affine space $\mathcal O_B$, hence irreducible.
For either braid word, the minimal-gallery endpoint map
identifies its open Bott--Samelson variety with
$\mathcal O_A$. The retraction strata partition
$\mathcal O_A\cap\mathcal O_B$ into locally closed products
$\mathbb A^{c^+}\times(\mathbb G_m)^{f^+}$, by the same
rank-one calculation as \cite[Lemma~5.4]{Hai24}.
Here $c^+$ and $f^+$ count the positive crossings and folds
of the block. Consequently, for either word there is exactly one
positively folded block of maximum dimension from $A$ to $B$.

It remains to prove independence of the sequence.
Suppose $\mathbf i$ and $\mathbf j$ differ by one braid move,
and put
$ 
C=\overline{S_{\mathbf i}(\xi)}
$ where the closure is taken in $Z(x,y)$. Since $\dim\xi=d$, the variety $C$ is an irreducible
component of $Z(x,y)$. The partition for $\mathbf j$ has a
unique stratum $S_{\mathbf j}(\eta)$ containing the generic
point of $C$. Its closure contains $C$ and has dimension
at most $d$, so
\[
\overline{S_{\mathbf j}(\eta)}=C,
\qquad \dim\eta=d.
\]
Both strata are dense open subsets of $C$.
For $z$ in their intersection, replace the rank-two
subgallery of $\widetilde\delta_{\mathbf i}(z)$ by the
minimal subgallery of the other braid type with the same
endpoints. The resulting gallery has type $\mathbf j$
and endpoint $z$, so it equals
$\widetilde\delta_{\mathbf j}(z)$ by uniqueness.
Their retractions $\xi$ and $\eta$ therefore agree outside
the braid block. Since their total dimensions agree, their
block dimensions agree.
The local uniqueness just proved identifies $\eta$ with
the braid replacement of $\xi$.

Every braid move therefore preserves $C$. For a fixed final
reduced expression, there is only one stratum containing
the generic point of $C$, so every sequence gives the
same final walk.
\end{proof}

\begin{rmk}
For explicit rank-two braid-move replacement algorithms, one may consult
\cite[Table~1 and the proof of Lemma~4.34]{MST19}.
\end{rmk}

\paragraph{The reduced expressions.}
Let $\de\in B_\mu(we,\lambda)$, let
$\ga=\operatorname{sk}(\de)$, and put $r=\ell(v)$.
Set $x_s=s_{t_1}\cdots s_{t_p}$, so that
$x_s\omega_\mu=(t_{\mu^*})_0$.
Lemma~\ref{lem:partial-conjugation} gives
\[
vx_s=x_s\mu_s(v),\qquad
\mu_s(v)\in W_{n_s},\qquad
\ell(\mu_s(v))=r.
\]
Choose reduced expressions
$v=a_1\cdots a_r$ and
$\mu_s(v)=a_1^{(s)}\cdots a_r^{(s)}$, and fix a reduced
expression for $w$.
Then
\[
w\,s_{t_1}\cdots s_{t_p}\,
a_1^{(s)}\cdots a_r^{(s)}\,\omega_\mu
\quad\text{and}\quad
w\,a_1\cdots a_r\,s_{t_1}\cdots s_{t_p}\,\omega_\mu
\]
are reduced expressions for the same element
$wv(t_{\mu^*})_0$.

Start from the last alcove of $\de$ and read the letters
$a_1^{(s)},\ldots,a_r^{(s)}$ in this order.
At step $j$, take the panel of the current alcove of type
$a_j^{(s)}$. Cross its wall if the current alcove and
$\De_U$ lie on the same side; otherwise fold at the panel.
Use the resulting last alcove for the next step.
Every appended step is a positive crossing or a positive
fold. Since the letters lie in $W_{n_s}$, all these steps
fix the terminal face, and the endpoint vertex remains
$-\lambda$.

Prepend the fixed minimal segment from $\Omega_0$ to
$w\Omega_0$, and denote the resulting alcove-to-alcove walk
by $\bar\de$. Its type is the first displayed reduced
expression. The prepended crossings replace the
initial-vertex contribution $\ell(w)$, while the terminal
segment contributes $r$. Hence
\[
\dim\bar\de=\dim\de+r
=\langle\rho,\lambda+\mu\rangle+\ell(v).
\]

\paragraph{Changing the reduced expression.}
Choose any sequence of braid moves from the first displayed
expression to the second. At every move, replace the
corresponding subwalk by a positively folded subwalk of the
new type with the same initial and terminal alcoves and the
same dimension, as in Lemma~\ref{lem:braid-existence}. The final walk has type
$w\,a_1\cdots a_r\,s_{t_1}\cdots s_{t_p}\,\omega_\mu$.
Remove its initial minimal segment from $\Omega_0$ to
$wv\Omega_0$, as in Lemma~\ref{lem:initial-minimal-segment},
and call the remaining vertex-to-vertex walk
$\de_{\mathrm{out}}$.

\begin{prop}\label{prop:braid-criterion}
The construction works for every sequence of braid moves
between the two displayed reduced expressions.
At every move the dimension-preserving positively folded
replacement is unique, and
\[
\de_{\mathrm{out}}=T^w_{v\leftarrow e}(\de).
\]
In particular, any two such sequences give the same output.
\end{prop}

\begin{proof}
Use the projection
of Section~4:
\[
q:Y_{wv}=Iwv(t_{\mu^*})_0I/I
\longrightarrow Y_w=Iw(t_{\mu^*})_0K/K.
\]
Its fibers are affine spaces of dimension $r=\ell(v)$,
and the root-group coordinates of Section~4 identify it
with a product projection.
Let $S_\de\subset Z_{\lambda,\mu}\cap Y_w$ be the
maximal-dimensional stratum indexed by $\de$, and let
$S_{\bar\de}\subset Y_{wv}$ be the stratum indexed by
the extended walk.

The terminal letters fix the terminal vertex.
Dropping them from a minimal lift of $\bar\de$ therefore
gives a lift of $\de$ with the same image in $Gr_G$.
Thus
\[
S_{\bar\de}\subset q^{-1}(S_\de).
\]
Both sides have dimension
$\langle\rho,\lambda+\mu\rangle+r$, and the right-hand
side is irreducible, being an affine-space bundle over
$S_\de$. Hence $S_{\bar\de}$ is dense in $q^{-1}(S_\de)$.

Write $y\Omega_0$ for the terminal alcove of $\bar\de$
and put $x=wv(t_{\mu^*})_0$.
Since this alcove has terminal vertex $-\lambda$,
\[
LUyI/I\cap IxI/I
\subset q^{-1}(Z_{\lambda,\mu}\cap Y_w).
\]
The right-hand side has dimension at most
$\langle\rho,\lambda+\mu\rangle+r$, by Section~4.
Since the left-hand side contains $S_{\bar\de}$ of this
dimension, $\bar\de$ satisfies the maximality hypothesis
of Lemma~\ref{lem:braid-existence}.
Every braid move is therefore defined, its replacement
is unique, and all sequences preserve the same irreducible
component of $LUyI/I\cap IxI/I$.

Let $S_{\mathrm{out}}$ be the stratum for the final
alcove-to-alcove walk, before removing its initial segment.
Component preservation and the preceding density give
\[
\overline{q(S_{\mathrm{out}})}^{\,Z_{\lambda,\mu}}
=
\overline{q(S_{\bar\de})}^{\,Z_{\lambda,\mu}}
=
\overline{S_\de}^{\,Z_{\lambda,\mu}}.
\]
Removing the initial segment gives
$\de_{\mathrm{out}}\in B_\mu(wv,\lambda)$.
The compatibility of the skeleton bijections with these
component labels, proved in
Proposition~\ref{prop:skeleton-bijection}, now yields
$\operatorname{sk}(\de_{\mathrm{out}})=\ga$.
By the definition of $T^w_{v\leftarrow e}$,
$\de_{\mathrm{out}}=T^w_{v\leftarrow e}(\de)$.
\end{proof}

\section{Raising and lowering operators on
\texorpdfstring{$B_\mu^v$}{the alcove-walk model}}\label{sec:crystal}

Fix $v\in W_{\mu^*}$. Proposition~\ref{prop:skeleton-bijection}
identifies $B_\mu^v$ with the LS galleries of type $\ga_\mu$.
We use this identification to give the set of maximal-dimensional
positively folded alcove walks a crystal structure.

We first transport the Gaussent--Littelmann operators through the
skeleton bijection. We then express the transported operators
as reflections and translations of alcove walks. The only
correction occurs when the direct operation changes the
stabilizer factor of the initial alcove.

We use the root operators of Section~3.2 on LS galleries. In
particular, a gallery ending at $-\lambda$ has weight $\lambda$.
Adjoin a symbol $\mathbf0$, distinct from every gallery or walk,
and extend every skeleton bijection by
$\operatorname{PR}_v(\mathbf0)=\mathbf0$.

\begin{dfn}\label{def:transported-operators}
For $\de\in B_\mu^v$, define
\begin{align}
e_\al(\de)
&:=\operatorname{PR}_v^{-1}
   \bigl(e_\al^{GL}(\operatorname{PR}_v(\de))\bigr),\label{eq:transport-e}\\
f_\al(\de)
&:=\operatorname{PR}_v^{-1}
   \bigl(f_\al^{GL}(\operatorname{PR}_v(\de))\bigr).\label{eq:transport-f}
\end{align}
An undefined gallery operator is interpreted as $\mathbf0$.
Set $e_\al(\mathbf0)=f_\al(\mathbf0)=\mathbf0$.
For $\de\in B_\mu^v(\lambda)$, put
\[
\operatorname{wt}(\de)=\lambda,\qquad
\varepsilon_\al(\de)
=\varepsilon_\al^{GL}(\operatorname{sk}(\de)),\qquad
\varphi_\al(\de)
=\varphi_\al^{GL}(\operatorname{sk}(\de)).
\]
\end{dfn}

\begin{thm}[Crystal structure on maximal-dimensional alcove walks]
\label{thm:alcove-crystal}
The preceding maps give $B_\mu^v$ the structure of the highest-weight
crystal $B(\mu)$ for $\widehat G$. The map
\[
\operatorname{PR}_v:B_\mu^v
\xrightarrow{\ \sim\ }\Gamma_{LS}^+(\ga_\mu)
\]
is an isomorphism of crystals. In particular,
\begin{align*}
\operatorname{sk}(e_\al\de)
&=e_\al^{GL}(\operatorname{sk}(\de)),\\
\operatorname{sk}(f_\al\de)
&=f_\al^{GL}(\operatorname{sk}(\de)).
\end{align*}
For nonzero values of the operators,
\begin{align*}
\operatorname{wt}(e_\al\de)
&=\operatorname{wt}(\de)+\al^\vee,&
\operatorname{wt}(f_\al\de)
&=\operatorname{wt}(\de)-\al^\vee,\\
\dim(e_\al\de)&=\dim\de+1,&
\dim(f_\al\de)&=\dim\de-1.
\end{align*}
Moreover,
\[
e_\al\de=\de'
\quad\Longleftrightarrow\quad
f_\al\de'=\de
\qquad(\de,\de'\in B_\mu^v),
\]
and
\[
\varphi_\al(\de)
=\varepsilon_\al(\de)+\langle\al,\operatorname{wt}(\de)\rangle.
\]
\end{thm}

\begin{proof}
The LS galleries, equipped with their root operators, form the
crystal $B(\mu)$ by \cite{GL05}, with the conventions of
Section~3.2. Proposition~\ref{prop:skeleton-bijection} is a
weight-preserving bijection. Transport of all the crystal maps
along this bijection proves the assertion and the displayed
operator relations. The dimension identities follow from
Definition~\ref{def:alcove-model} and $\langle\rho,\al^\vee\rangle=1$.
\end{proof}

The crystal structure is thus obtained from the
Gaussent--Littelmann theorem. The remaining task is to describe
its operators directly on alcove walks.

For later use, let $\ga=\operatorname{sk}(\de)$ and use the
height $m_\al(\ga)$ of Definition~\ref{def:gallery-root-operators}.
The gallery formulas give
\begin{equation}\label{eq:gallery-string-lengths}
\varepsilon_\al(\de)=m_\al(\ga),\qquad
\varphi_\al(\de)=m_\al(\ga)+\langle\al,\lambda\rangle
\quad(\de\in B_\mu^v(\lambda)).
\end{equation}
Thus $e_\al\de=\mathbf0$ exactly when $m_\al(\ga)=0$, and
$f_\al\de=\mathbf0$ exactly when
$m_\al(\ga)=-\langle\al,\lambda\rangle$.
Replacing $m_\al(\ga)$ here by the corresponding maximum on
the alcove walk requires the height comparison discussed below.

\begin{cor}\label{cor:transition-crystal}
The bijection
\[
T_{v\leftarrow v'}
:=\operatorname{PR}_v^{-1}\circ\operatorname{PR}_{v'}
:B_\mu^{v'}\xrightarrow{\ \sim\ }B_\mu^v
\]
is an isomorphism of crystals. Its restriction to a fixed initial
direction is $T^w_{v\leftarrow v'}$ of
Proposition~\ref{prop:transition}.
\end{cor}
\begin{proof}
Both skeleton maps are crystal isomorphisms by
Theorem~\ref{thm:alcove-crystal}.
\end{proof}

\subsection{Bijection between the model of Haines and LS galleries}
\label{subsec:haines-inverse}

By Remark~\ref{rmk:haines-alcove-model}, $B_\mu^e$ is the vertex-to-vertex version of
Haines' maximal-dimensional alcove-walk model
\[
\bigsqcup_\lambda\mathcal M_\mu^{\De_U}(\lambda).
\]
As before, $\lambda$ ranges over the cocharacters in the convex
hull of $W_0\mu$ with $\mu-\lambda\in\mathbb Z\Phi^\vee$.
The map to LS galleries sends a walk to its skeleton.

The inverse can be described using a sequence of raising operators.
The fixed gallery $\ga_\mu$, which ends at $\mu^*$, is the
lowest-weight gallery, of weight $-\mu^*=w_0\mu$.
For $v\in W_{\mu^*}$, the minimal walk
\[
\de_{\mathrm{low}}^v:=v\de_\mu
\]
has skeleton $\ga_\mu$ and dimension $\ell(v)$; it is the
lowest-weight element of $B_\mu^v$.

\begin{prop}\label{prop:inverse-by-raising}
Suppose
\[
\ga=e_{\al_N}^{GL}\cdots e_{\al_1}^{GL}(\ga_\mu),
\]
where every intermediate gallery is defined. Then
\[
\operatorname{PR}_v^{-1}(\ga)
=e_{\al_N}\cdots e_{\al_1}(\de_{\mathrm{low}}^v).
\]
The result is independent of the chosen sequence of simple roots.
\end{prop}
\begin{proof}
Every LS gallery is obtained from the lowest-weight gallery by
a sequence of raising operators. Apply
Theorem~\ref{thm:alcove-crystal} at each step. The final walk has
skeleton $\ga$, so its uniqueness follows from
Proposition~\ref{prop:skeleton-bijection}.
\end{proof}

The explicit formulas in Section~\ref{subsec:explicit-formulas}
make this inverse procedure combinatorial. One may maintain the lifts
for every $v\in W_{\mu^*}$, using the transition maps when an exceptional
raising step occurs. Alternatively, start with
$\de_{\mathrm{low}}^{w_{0,\mu^*}}$, apply only the unmodified raising
formulas, and then apply $T^w_{e\leftarrow w_{0,\mu^*}}$, where
$w=t(\ga)$. Corollary~\ref{rmk:longest-stabilizer} justifies the latter
procedure.

\subsection{Direct GL operators and the skeleton}
\label{subsec:raw-operators}

For alcove walks, we write
\[
\mathsf E_\al\de=e_\al^{GL}(\de),\qquad
\mathsf F_\al\de=f_\al^{GL}(\de)
\]
for the operators of Definition~\ref{def:gallery-root-operators}.
Their domains, the height $m_\al(\de)$, and the cut indices $j,k$
are those of \eqref{eq:raw-raising} and \eqref{eq:raw-lowering},
applied to the faces $\De_i',\De_i$ of the walk.
The symbols $e_\al$ and $f_\al$ denote the transported crystal
operators of Definition~\ref{def:transported-operators}.
We now compare the direct operators with the skeleton map and
identify when they agree with the transported operators.

\begin{lem}\label{lem:skeleton-height}
Let $\de\in B_\mu^v(\lambda)$ and let
$\ga=\operatorname{sk}(\de)$. Then, for every simple root $\al$,
\[
m_\al(\de)=m_\al(\ga).
\]
\end{lem}

\begin{proof}
We work over an algebraic closure of $k$.
The proof of Proposition~\ref{prop:skeleton-bijection} gives
\[
\pi\bigl(C_{X_1}(\de)\bigr)\subset C_{X_2}(\ga),
\qquad
\overline{\pi\bigl(C_{X_1}(\de)\bigr)}
=C_{X_2}(\ga),
\]
where the closure is taken in $C_{X_2}(\ga)$.

Let $M$ be the standard Levi subgroup with root system
$\{\al,-\al\}$, and put
$ 
K_M=M(k[[t]]),A=\langle\al,\lambda\rangle.
$  
The quotient $U\to U\cap M=U_\al$ defines a homomorphism
\[
\chi_\al:LU\longrightarrow k((t)),
\qquad
\chi_\al(u_\beta(z))=
\begin{cases}
z,&\beta=\al,\\
0,&\beta\in\Phi^+\setminus\{\al\}.
\end{cases}
\]
Here $\al$ is simple. The rank-one (non-continuous) retraction $r_\al:Gr_G
\longrightarrow Gr_M$ satisfies
$ 
r_\al(ut^\lambda K)
=u_\al\bigl(\chi_\al(u)\bigr)t^\lambda K_M.
$ 

We recall the root-coordinate calculation. For $\eta=\de$ or $\ga$, write its faces as
$F_i'\subset F_i$, and set
\[
J_i(\eta)=
\left\{
(\beta,n)\in\Phi^+\times\mathbb Z:
F_i'\subset H_{\beta,n},\
F_i\subset H_{\beta,<n}
\right\}.
\]
Using the root-group coordinates of \cite[Section 10]{GL05}, with the retraction direction reversed,
the large faces of a point of the retraction stratum are
\begin{equation}\label{eq:prefix-root-coordinates}
\widetilde F_d=
\left(
\prod_{i=0}^{d}
\prod_{(\beta,n)\in J_i(\eta)}
u_\beta(z_{i,\beta,n}t^{-n})
\right)F_d.
\end{equation}
The products are taken in fixed root orders, and the
$z_{i,\beta,n}$ are independent affine coordinates.
Consequently,
\begin{equation}\label{eq:prefix-simple-root-coordinate}
\chi_\al\left(
\prod_{i=0}^{d}
\prod_{(\beta,n)\in J_i(\eta)}
u_\beta(z_{i,\beta,n}t^{-n})
\right)
=
\sum_{\substack{0\leq i\leq d\\
                 (\al,n)\in J_i(\eta)}}
z_{i,\al,n}t^{-n}.
\end{equation}

Put $m=m_\al(\eta)$. Since the terminal vertex is $-\lambda$,
$ 
m\geq-A.
$ 
For every integer $n$ with $-A<n\leq m$, the portion of
$\eta$ from a small face in $H_{\al,m}$ to $-\lambda$
contains a pair
$ 
F_i'\subset H_{\al,n},
\qquad
F_i\subset H_{\al,<n}.
$ 
Thus
\[
\{-A+1,\ldots,m\}
\subset
\{n:(\al,n)\in J_i(\eta)\text{ for some }i\}
\subset
\{n:n\leq m\}.
\]
Moreover,
\[
u_\al(z)t^\lambda K_M=t^\lambda K_M
\quad\Longleftrightarrow\quad
z\in t^A k[[t]].
\]
Applying \eqref{eq:prefix-simple-root-coordinate} at the
endpoint therefore gives
\[
\begin{aligned}
r_\al\bigl(f(C_{X_1}(\de))\bigr)
&=
\left\{
u_\al\left(
\sum_{n=-A+1}^{m_\al(\de)}c_nt^{-n}
\right)t^\lambda K_M:
c_n\in k
\right\},\\
r_\al\bigl(g(C_{X_2}(\ga))\bigr)
&=
\left\{
u_\al\left(
\sum_{n=-A+1}^{m_\al(\ga)}c_nt^{-n}
\right)t^\lambda K_M:
c_n\in k
\right\}.
\end{aligned}
\]
An empty sum gives the point $t^\lambda K_M$.
In particular, the dimensions of these two images are
$ 
A+m_\al(\de),A+m_\al(\ga),
$ 
respectively. This is the rank-one calculation of
\cite[Section 5.3]{BG07}, in our convention.

All endpoint images lie in $LUt^\lambda K/K$, on which
$r_\al$ is a morphism. Since $f=g\circ\pi$ and
$\pi(C_{X_1}(\de))$ is dense in $C_{X_2}(\ga)$,
$ 
\overline{r_\al\bigl(f(C_{X_1}(\de))\bigr)}
=
\overline{r_\al\bigl(g(C_{X_2}(\ga))\bigr)}.
$ 
Comparing dimensions gives
$ 
A+m_\al(\de)=A+m_\al(\ga),
$ 
as required.
\end{proof}

\begin{lem}\label{lem:walk-dimension-bound}
Let $\xi$ be a positively folded walk of type $t_{\de_\mu}$,
with first alcove $wv\Omega_0$ and endpoint $-\lambda$, where
$w\in W^{\mu^*}$ and $v\in W_{\mu^*}$. Then
\begin{equation}\label{eq:walk-dimension-bound}
\dim\xi\leq\langle\rho,\lambda+\mu\rangle+\ell(v).
\end{equation}
\end{lem}
\begin{proof}
Prepend the minimal segment from $\Omega_0$ to $wv\Omega_0$.
By \cite[Lemma 5.4]{Hai24}, the resulting walk indexes a nonempty
locally closed stratum of dimension $\dim\xi$ in
\[
\bigsqcup_{u\in W_0}
\bigl(LUt^\lambda uI/I\cap Iwv(t_{\mu^*})_0I/I\bigr).
\]
By \eqref{eq:fw-preimage} and the affine-space projection in
Section~4.2, this intersection is isomorphic to
$ 
\bigl(LUt^\lambda K/K\cap Iw(t_{\mu^*})_0K/K\bigr)
\times\mathbb A^{\ell(v)}.
$ 
The first factor is contained in $Z_{\lambda,\mu}$, whose
dimension is at most $\langle\rho,\lambda+\mu\rangle$;
see \cite[Proposition 3.2]{Hai24} with $P=B$.
This proves \eqref{eq:walk-dimension-bound}.
\end{proof}

\begin{lem}\label{lem:maximal-wall}
Let $\al$ be simple.
\begin{enumerate}
\item If $\ga$ is an LS gallery and $m=m_\al(\ga)$, then
$\Ga_i\subset H_{\al,\leq m}$ for every $i$.
\item If $\de\in B_\mu^v$ and $m=m_\al(\de)$, no crossing
at an index $a>0$ goes from $H_{\al,<m}$ to $H_{\al,>m}$.
In particular, if $m>0$, or if
$\De_0\subset H_{\al,<m}$, then
$\De_i\subset H_{\al,<m}$ for every $i$.
\end{enumerate}
\end{lem}
\begin{proof}
For (1), suppose that $a$ is the first index with
$\Ga_a\subset H_{\al,>m}$.
Then $\Ga_a'\subset H_{\al,m}$. Let $b>a$ be the next
index with $\Ga_b'\subset H_{\al,m}$, allowing $b=s+1$.
Such an index exists because the endpoint has integral
$\al$-coordinate at most $m$. Define
\[
\eta_i=
\begin{cases}
s_{\al,m}\Ga_i,&a\leq i<b,\\
\Ga_i,&i<a\text{ or }i\geq b.
\end{cases}
\]
The excursion calculation of \cite[Section 6]{GL05}, conjugated
by $x\mapsto-x$ as in Section~3.2, gives
$\eta\in\Ga^+(\lambda,\ga_\mu)$ and
$\dim\eta=\dim\ga+1$. This contradicts the gallery dimension
bound $\dim\eta\leq\langle\rho,\lambda+\mu\rangle=\dim\ga$.

For (2), suppose such a crossing occurs at $a>0$, and let
$b>a$ be the next index with $\De_b'\subset H_{\al,m}$,
allowing $b=p+1$. Since no small face lies in $H_{\al,m+1}$,
\[
\De_i\subset\{x:m<\langle\al,x\rangle<m+1\}
\qquad(a\leq i<b).
\]
Replace $\De_i$ by $s_{\al,m}\De_i$ for $a\leq i<b$ and
leave the other alcoves unchanged. There is no interior
$\al$-panel. Every other positive fold remains positive, since
$s_\al(\Phi^+\setminus\{\al\})=\Phi^+\setminus\{\al\}$.
At $a$ a negative crossing becomes a positive fold, increasing
the dimension by one. If $b\leq p$, the positive crossing at
$b$ becomes a positive fold, with no change in its contribution;
if $b=p+1$, the terminal vertex is fixed. The modified walk
has the same first alcove and endpoint and dimension
$\dim\de+1$, contrary to Lemma~\ref{lem:walk-dimension-bound}.
Finally, $m>0$ implies $\De_0\subset H_{\al,<m}$.
\end{proof}

\begin{prop}\label{prop:cut-positions}
Let $\de\in B_\mu^v$ and $\ga=\operatorname{sk}(\de)$.
If $m_\al(\de)>0$, then
\[
\operatorname{sk}(\mathsf E_\al\de)=e_\al^{GL}\ga,
\]
whether or not $\mathsf E_\al\de$ belongs to $B_\mu^v$.
If $\mathsf F_\al\de\in B_\mu^v$, then
$\operatorname{sk}(\mathsf F_\al\de)=f_\al^{GL}\ga$.
\end{prop}

\begin{proof}
We apply Definition~\ref{def:gallery-root-operators} to both
$\de$ and $\ga$. Thus $\mathsf E_\al\de=e_\al^{GL}(\de)$ and
$\mathsf F_\al\de=f_\al^{GL}(\de)$. Since these operators preserve
type, the resulting alcove walks still have type $t_{\de_\mu}$,
so their skeletons are defined without requiring their values to
belong to $B_\mu^v$.

We first prove the raising assertion. Write
$\ga:\Ga_0'\subset\Ga_0\supset\cdots\supset\Ga_{s+1}'$.
By Lemma~\ref{lem:skeleton-height}, put
$m=m_\al(\de)=m_\al(\ga)>0$ and $h=m-1$.
Let $j<k$ be the raising cut indices for $\de$ in
Definition~\ref{def:gallery-root-operators}, and let $l<q$ be
those for $\ga$; explicitly,
\[
q=\min\{i:\Ga_i'\subset H_{\al,m}\},\qquad
l=\max\{i<q:\Ga_i'\subset H_{\al,h}\}.
\]
Let $a,b$ be determined by
$i_{a-1}<j\leq i_a$ and $i_{b-1}<k\leq i_b$,
with $a=0$ if $j=0$ and $b=s+1$ if $k=p+1$.
The paragraph following Definition~\ref{def:skeleton} gives
$\Ga_a'\subset\overline{\De_j'}\subset H_{\al,h}$ and
$\Ga_b'\subset\overline{\De_k'}\subset H_{\al,m}$.
Thus $a<b$ and $q\leq b$.

We first prove
\begin{equation}\label{eq:upper-cut-wall}
\Ga_i\subset H_{\al,m}\qquad(q\leq i<b).
\end{equation}
Since $k$ is the first index in $H_{\al,m}$ and $m>0$,
$\De_t\subset H_{\al,<m}$ for $t<k$; hence
$\Ga_i\subset H_{\al,\leq m}$ for $i<b$.
If \eqref{eq:upper-cut-wall} fails, let
$d=\min\{q\leq i<b:\Ga_i\not\subset H_{\al,m}\}$.
Then $\Ga_d'\subset H_{\al,m}$ and
$\Ga_d\subset H_{\al,<m}$.

For the geometric argument, work over an algebraic closure of $k$
and use $\chi_\al$ from the proof of Lemma~\ref{lem:skeleton-height}.
The root-group filtration gives
\[
\chi_\al\bigl(LU\cap\operatorname{Stab}(\Ga_d)\bigr)
=t^{\left\lceil-\inf_{x\in\Ga_d}\langle\al,x\rangle\right\rceil}k[[t]]
\subset t^{1-m}k[[t]].
\]
Consequently
$c_d(u\Ga_d):=[t^{-m}]\chi_\al(u)$ is well defined on
$LU\cdot\Ga_d$ and is regular on its finite-dimensional
root-coordinate charts. Evaluation at the $d$th large face satisfies
\[
\operatorname{ev}_d(h_0,\ldots,h_s)=h_0\cdots h_dF_{n_d},
\qquad
\operatorname{ev}_d\pi(g_0,\ldots,g_p)
=g_0\cdots g_{i_d}F_{n_d}.
\]
Since $i_d<k$, this walk prefix has no panel in $H_{\al,m}$.
Its initial root coordinates have exponent $0$, and $m>0$.
Thus \eqref{eq:prefix-simple-root-coordinate} gives
$c_d\circ\operatorname{ev}_d\circ\pi=0$ on $C_{X_1}(\de)$.
On $C_{X_2}(\ga)$, however, the same formula gives
\[
c_d\circ\operatorname{ev}_d
=\sum_{\substack{0\leq i\leq d\\(\al,m)\in J_i(\ga)}}
z_{i,\al,m}.
\]
This function is nonzero, since $(\al,m)\in J_d(\ga)$ and
$z_{d,\al,m}$ is an independent coordinate. This contradicts
the density of $\pi(C_{X_1}(\de))$ in $C_{X_2}(\ga)$.
Therefore \eqref{eq:upper-cut-wall} holds.

It follows that $a<q$: otherwise \eqref{eq:upper-cut-wall}
would give $\Ga_a'\subset H_{\al,m}\cap H_{\al,h}$.
Hence $a\leq l<q\leq b$.
The definitions of $j,k$ also give
\[
\Ga_i\subset\{x:h\leq\langle\al,x\rangle\leq m\}
\qquad(a\leq i<b).
\]
Put $\ga^+=e_\al^{GL}\ga$. This is an LS gallery with
$m_\al(\ga^+)=h$, so Lemma~\ref{lem:maximal-wall}(1) gives
$\ga_i^+\subset H_{\al,\leq h}$ for all $i$.
For $i<l$, the gallery raising formula gives $\ga_i^+=\Ga_i$.
Consequently
\begin{equation}\label{eq:lower-cut-wall}
\Ga_i\subset H_{\al,h}\qquad(a\leq i<l).
\end{equation}

By Definition~\ref{def:skeleton}, the large faces of
$\operatorname{sk}(\mathsf E_\al\de)$ are
\begin{equation}\label{eq:induced-skeleton-raising}
\bigl(\operatorname{sk}(\mathsf E_\al\de)\bigr)_i=
\begin{cases}
\Ga_i,&i<a,\\
s_{\al,h}\Ga_i,&a\leq i<b,\\
\Ga_i-\al^\vee,&i\geq b.
\end{cases}
\end{equation}
On $H_{\al,h}$ the reflection $s_{\al,h}$ is the identity;
on $H_{\al,m}$ it agrees with translation by $-\al^\vee$.
Equations \eqref{eq:upper-cut-wall} and \eqref{eq:lower-cut-wall}
therefore allow $a,b$ in \eqref{eq:induced-skeleton-raising}
to be replaced by $l,q$. By \eqref{eq:raw-raising}, these are
exactly the large faces of $e_\al^{GL}(\ga)$.

For the smaller skeleton faces, the induced maps are the identity
for $i\leq a$, reflection in $H_{\al,h}$ for $a<i\leq b$,
and translation by $-\al^\vee$ for $i>b$.
The same wall containments, together with incidence, give
\[
\Ga_i'\subset H_{\al,h}\quad(a\leq i\leq l),\qquad
\Ga_i'\subset H_{\al,m}\quad(q\leq i\leq b).
\]
The maps therefore also agree on the smaller faces when the cuts
$a,b$ are replaced by $l,q$. This includes the boundary cases:
if $j=0$, then $m=1$ and the reflection fixes the initial vertex;
if $k=p+1$, reflection of the terminal vertex in $H_{\al,m}$
across $H_{\al,h}$ agrees with translation by $-\al^\vee$.
Both galleries end at $-(\lambda+\al^\vee)$. Hence
\[
\operatorname{sk}(\mathsf E_\al\de)
=\operatorname{sk}(e_\al^{GL}(\de))
=e_\al^{GL}(\ga).
\]

For the lowering assertion, put
$\xi=\mathsf F_\al\de=f_\al^{GL}(\de)\in B_\mu^v$.
Write $M=m_\al(\de)$ and let $K<J$ be the lowering cut indices
for $\de$ in Definition~\ref{def:gallery-root-operators}.
The same definition gives $m_\al(\xi)=M+1$: the first small
face of $\xi$ in $H_{\al,M+1}$ has index $J$, and the last
preceding small face in $H_{\al,M}$ has index $K$.
Thus the raising operation on $\xi$ reverses the reflection
and translation in the lowering operation on $\de$, giving
\[
e_\al^{GL}(\xi)=\de.
\]
Since $M+1>0$ and $\xi\in B_\mu^v$, the raising assertion applies:
\[
e_\al^{GL}(\operatorname{sk}(\xi))
=\operatorname{sk}(e_\al^{GL}(\xi))
=\ga.
\]
Both skeletons are LS galleries by
Proposition~\ref{prop:skeleton-bijection}. Applying the inverse
operator $f_\al^{GL}$ on LS galleries yields
\[
\operatorname{sk}(\mathsf F_\al\de)
=\operatorname{sk}(\xi)=f_\al^{GL}(\ga),
\]
as required.
\end{proof}

\begin{prop}[Direct operator comparison]\label{prop:direct-reduction}
If $\mathsf E_\alpha\de$ belongs to $B_\mu^v$, then
$\mathsf E_\alpha\de=e_\alpha\de$.
If $\mathsf F_\alpha\de$ belongs to $B_\mu^v$, then
$\mathsf F_\alpha\de=f_\alpha\de$.
\end{prop}
\begin{proof}
The raising assertion follows from Proposition~\ref{prop:cut-positions}
and uniqueness of the maximal lift.
For lowering, put $\xi=\mathsf F_\alpha\de$.
The reflection-and-translation formulas give
$\mathsf E_\alpha\xi=\de$; these are the formal inverse identities
of \cite[Section 6]{GL05}, in our convention.
Apply the raising assertion to $\xi$ and then apply $f_\alpha$.
\end{proof}

\subsection{Explicit folding formulas}\label{subsec:explicit-formulas}

Let
$ 
\xi:\Xi_0'\subset\Xi_0\supset\Xi_1'\subset\cdots
\subset\Xi_p\supset\Xi_{p+1}'
$ 
be an alcove walk. If $\Xi_d'\subset H_{\alpha,0}$, let
$R_{\alpha,d}\xi$ be the reflection of its prefix:
\begin{equation}\label{eq:prefix-reflection}
(R_{\alpha,d}\xi)_i=
\begin{cases}
s_\alpha\Xi_i,&i<d,\\
\Xi_i,&i\geq d.
\end{cases}
\end{equation}
The joining face is fixed by $s_\alpha$.

\begin{prop}[Raising]\label{prop:raising-criterion}
Let $\de\in B_\mu(wv,\lambda)$, let $\ga=\operatorname{sk}(\de)$,
and use the raising indices $m,j,k$ of \eqref{eq:raw-raising}.
The operator is undefined when $m=0$.
If $m>0$, the only exceptional case is
\[
m=1,\qquad j=0,\qquad (s_\alpha wv)^{\mu^*}=w.
\]
In this case write $s_\alpha w=ws_\beta$, with $s_\beta\in W_{\mu^*}$ simple,
and put $v^+=s_\beta v$. Then $\ell(v^+)=\ell(v)+1$.
Set
\[
\widehat\de=T^w_{v^+\leftarrow v}(\de),\qquad
\widetilde\de=\mathsf E_\alpha\widehat\de,
\]
and let $d$ be the first positive fold of $\widetilde\de$ in
$H_{\alpha,0}$. Such a fold exists. The formula is
\begin{equation}\label{eq:candidate-raising}
e_\alpha\de=
\begin{cases}
R_{\alpha,d}\bigl(\mathsf E_\alpha
T^w_{v^+\leftarrow v}(\de)\bigr),&\text{in the exceptional case},\\
\mathsf E_\alpha\de,&\text{otherwise}.
\end{cases}
\end{equation}
\end{prop}
\begin{proof}
Suppose $m>0$. By Lemma~\ref{lem:maximal-wall}(2),
$\De_i\subset H_{\alpha,<m}$ for every $i$.
Thus the step at $k$ is a positive fold, unless $k=p+1$.
For $j\leq i<k$ one has
$m-1<\langle\alpha,\De_i\rangle<m$.
There is no interior $\alpha$-panel, and
$s_\alpha(\Phi^+\setminus\{\alpha\})=\Phi^+\setminus\{\alpha\}$
preserves every other folding sign. If $j>0$, the negative
crossing at $j$ becomes a positive fold. If $k\leq p$, the
positive fold at $k$ becomes a positive crossing.
The changes in their dimension contributions are respectively
$+1$ and $0$. If $j=0$, then $m=1$ and
$\ell(s_\alpha wv)=\ell(wv)+1$, so the initial contribution
increases by one instead. Hence $\mathsf E_\alpha\de$ is
positively folded and
$\dim(\mathsf E_\alpha\de)=\dim\de+1$.

If $j>0$, the first alcove is unchanged. If $j=0$,
Lemma~\ref{lem:coset-alternative} gives either stabilizer factor
$v$, or $s_\alpha w=ws_\beta$ with
$v^+=s_\beta v$ and $\ell(v^+)=\ell(v)+1$.
In the first case $\mathsf E_\alpha\de\in
B_\mu^v(\lambda+\alpha^\vee)$, and
Proposition~\ref{prop:direct-reduction} proves the formula.

Now assume the exceptional case. The walk
$\widehat\de=T^w_{v^+\leftarrow v}(\de)$ has skeleton $\ga$,
height one, and first alcove $s_\alpha wv\Omega_0\subset
H_{\alpha,<0}$. Denote its raising indices by $\widehat j,\widehat k$.
Then $\widehat j>0$, and the step at $\widehat j$ is a
negative crossing of $H_{\alpha,0}$. The preceding case gives
\[
\widetilde\de=\mathsf E_\alpha\widehat\de
\in B_\mu(wv^+,\lambda+\alpha^\vee),
\qquad
\operatorname{sk}(\widetilde\de)=e_\alpha^{GL}\ga.
\]
The step at $\widehat j$ becomes a positive fold in $H_{\alpha,0}$.
Thus the first such fold $d$ exists and $d\leq\widehat j$.

Use the gallery index $l$ from the proof of
Proposition~\ref{prop:cut-positions}. For the original walk,
$j=0$, so \eqref{eq:lower-cut-wall} gives
$\Ga_i\subset H_{\alpha,0}$ for $0\leq i<l$;
also $\Ga_l'\subset H_{\alpha,0}$.
For $\widehat\de$, let $\widehat a$ be determined by
$i_{\widehat a-1}<\widehat j\leq i_{\widehat a}$.
The same proposition gives $\widehat a\leq l$.
Since $d\leq\widehat j\leq i_l$, every large skeleton face
changed by $R_{\alpha,d}$ has index less than $l$.
These faces are fixed by $s_\alpha$, and therefore
\[
\operatorname{sk}(R_{\alpha,d}\widetilde\de)
=\operatorname{sk}(\widetilde\de)=e_\alpha^{GL}\ga.
\]

By Lemma~\ref{lem:skeleton-height},
$m_\alpha(\widetilde\de)=0$. Its first alcove lies in
$H_{\alpha,<0}$, so Lemma~\ref{lem:maximal-wall}(2) gives
$\widetilde\De_i\subset H_{\alpha,<0}$ for all $i$.
The gallery formula gives $(e_\alpha^{GL}\ga)'_r=\Ga_r'$
for $r\leq l$, since $\Ga_l'\subset H_{\alpha,0}$.
For $i<d$, the paragraph following Definition~\ref{def:skeleton} gives
$\Ga_r'\subset\overline{\widetilde\De_i}$ for some $r\leq l$;
these small faces lie in $H_{\alpha,0}$. Hence
$-1<\langle\alpha,\widetilde\De_i\rangle<0$ for $i<d$.
There is no $\alpha$-crossing before $d$, no positive fold in
$H_{\alpha,-1}$, and no fold in $H_{\alpha,0}$ by the choice
of $d$. Prefix reflection therefore preserves all folding
signs before $d$. At $d$ it replaces a positive fold by a
positive crossing, with unchanged contribution.
The initial contribution decreases by one. Thus
\[
\dim R_{\alpha,d}\widetilde\de
=\dim\widetilde\de-1
=\langle\rho,\lambda+\mu+\alpha^\vee\rangle+\ell(v).
\]
Its first alcove is $wv\Omega_0$, its endpoint is
$-(\lambda+\alpha^\vee)$, and its skeleton is $e_\alpha^{GL}\ga$.
Proposition~\ref{prop:skeleton-bijection} proves the assertion.
\end{proof}

The lowering operator has the following analogous formula.

\begin{prop}[Lowering]\label{prop:lowering-criterion}
Let $\de\in B_\mu(wv,\lambda)$ and
$\ga=\operatorname{sk}(\de)$. Use the lowering indices $m,k,j$
of \eqref{eq:raw-lowering}.
The operator is undefined when $m=-\langle\alpha,\lambda\rangle$.
Otherwise the only exceptional case is a positive crossing at
$\De_k'\subset H_{\alpha,m}$.
Then
\[
m=0,\qquad s_\alpha w=ws_\beta,\qquad
v^+=s_\beta v,\qquad \ell(v^+)=\ell(v)+1
\]
for a simple reflection $s_\beta\in W_{\mu^*}$.
The prefix-reflected walk $\de'=R_{\alpha,k}\de$
belongs to $B_\mu(wv^+,\lambda)$ and has skeleton $\ga$.
The formula is
\begin{equation}\label{eq:candidate-lowering}
f_\alpha\de=
\begin{cases}
T^w_{v\leftarrow v^+}\bigl(\mathsf F_\alpha\de'\bigr),
 &\text{if the last visit to }H_{\alpha,m}\text{ is a positive crossing},\\
\mathsf F_\alpha\de,&\text{otherwise}.
\end{cases}
\end{equation}
\end{prop}

\begin{proof}
The assertion about the domain follows from
Lemma~\ref{lem:skeleton-height} and \eqref{eq:gallery-string-lengths}.
Suppose the operator is defined. By the choice of $k,j$,
\[
m-1<\langle\alpha,\De_i\rangle<m\qquad(k\leq i<j).
\]
Thus, if the step at $k>0$ is a crossing, it is positive.
Assume first that $k=0$ or that the step at $k$ is a positive
fold. If $k>0$, reflection replaces this fold by a negative
crossing, decreasing its contribution by one. If $k=0$, it
decreases the initial contribution by one, since
$\De_0\subset H_{\alpha,<0}$ and
$\ell(s_\alpha wv)=\ell(wv)-1$.
If $j\leq p$, the positive crossing at $j$ becomes a positive
fold; its contribution remains one. There is no interior
$\alpha$-panel, and all other folding signs are preserved.
Consequently $\mathsf F_\alpha\de$ is positively folded and
$\dim(\mathsf F_\alpha\de)=\dim\de-1$.

For $k>0$ its first alcove is unchanged. For $k=0$, the
finite parabolic coset decomposition gives either the same
stabilizer factor $v$, or a factor $v^-$ with
$\ell(v^-)=\ell(v)-1$. The latter is impossible by
Lemma~\ref{lem:walk-dimension-bound}, since it would give
\[
\dim(\mathsf F_\alpha\de)
=\langle\rho,\lambda+\mu-\alpha^\vee\rangle+\ell(v)
>\langle\rho,\lambda+\mu-\alpha^\vee\rangle+\ell(v^-).
\]
Thus $\mathsf F_\alpha\de\in B_\mu^v(\lambda-\alpha^\vee)$,
and Proposition~\ref{prop:direct-reduction} applies.

It remains to consider a positive crossing at $k$.
Lemma~\ref{lem:maximal-wall}(2) forces $m=0$ and
\[
0<\langle\alpha,\De_i\rangle<1\qquad(0\leq i<k);
\]
otherwise the portion above $H_{\alpha,m}$ would begin with
a crossing at a positive index. In particular, $k$ is the
first panel index in $H_{\alpha,0}$.
Put $\de'=R_{\alpha,k}\de$, and let $b$ satisfy
$i_{b-1}<k\leq i_b$.
By Lemma~\ref{lem:maximal-wall}(1) and
Lemma~\ref{lem:skeleton-height}, $\Ga_i\subset H_{\alpha,\leq0}$.
For $i<b$, the face $\Ga_i$ is contained in an alcove of
the prefix, so also $\Ga_i\subset H_{\alpha,\geq0}$.
Therefore $\Ga_i\subset H_{\alpha,0}$ for $i<b$,
and $\Ga_b'\subset\De_k'\subset H_{\alpha,0}$.
Definition~\ref{def:skeleton} now gives
$\operatorname{sk}(\de')=\ga$.

The crossing at $k$ becomes a positive fold.
There is no other $\alpha$-panel in the reflected prefix,
and its non-$\alpha$ folding signs are preserved.
Since $\ell(s_\alpha wv)=\ell(wv)+1$, we obtain
$\dim\de'=\dim\de+1$.
Moreover $s_\alpha\Ga_0=\Ga_0=w\Upsilon_0$, so
$(s_\alpha wv)^{\mu^*}=w$.
The parabolic coset decomposition therefore gives
$s_\alpha w=ws_\beta$ with $s_\beta\in W_{\mu^*}$ simple and
$\ell(s_\beta v)=\ell(v)+1$. Hence
$\de'\in B_\mu(wv^+,\lambda)$, where $v^+=s_\beta v$.

The last visit of $\de'$ to $H_{\alpha,0}$ is the positive
fold at $k$, so the ordinary case just proved gives
\[
\mathsf F_\alpha\de'\in B_\mu(wv^+,\lambda-\alpha^\vee),
\qquad
\operatorname{sk}(\mathsf F_\alpha\de')=f_\alpha^{GL}\ga.
\]
Its first alcove remains $wv^+\Omega_0$.
Applying $T^w_{v\leftarrow v^+}$ preserves this skeleton and
gives the required element of $B_\mu(wv,\lambda-\alpha^\vee)$.
\end{proof}

\begin{cor}\label{rmk:longest-stabilizer}
For $v=w_{0,\mu^*}$, neither exceptional case occurs. Hence
\[
e_\alpha\de=\mathsf E_\alpha\de,\qquad
f_\alpha\de=\mathsf F_\alpha\de
\quad(\de\in B_\mu^{w_{0,\mu^*}})
\]
on their respective domains.
\end{cor}
\begin{proof}
Both exceptional cases require an element of $W_{\mu^*}$ with length
$\ell(w_{0,\mu^*})+1$, which is impossible. This is precisely the alcove walks used in \cite{BG07}, adapted to our notation.
\end{proof}

Together with Proposition~\ref{prop:inverse-by-raising}, these formulas
give the two inverse procedures described in
Section~\ref{subsec:haines-inverse}.

\subsection{Duality}\label{subsec:duality}

Let $\de$ be a walk from $0$ to $-\lambda$. Translate every face
of $\de$ by $+\lambda$, and then read the translated face sequence
in the opposite order: start with the translated terminal vertex,
then the translated last alcove, its preceding panel, and so on,
finishing at the translated initial vertex. Repeated alcoves are
retained, so folds remain folds. Denote the resulting walk by
$\de^*$. It starts at $0$ and ends at $\lambda$.

To choose the dual minimal walk, first apply $w_{0,\mu^*}$ to
every face of $\de_\mu$ and translate every face by $-\mu^*$.
Read this face sequence in the opposite order as above, and finally
apply $w_0$ to every face. Call the resulting walk $\de_{\mu^*}$.
It starts at $0$, has first alcove $\Omega_0$, and ends at $\mu$:
the terminal alcove of $w_{0,\mu^*}\de_\mu$ is
$t_{\mu^*}w_0\Omega_0$, which becomes $w_0\Omega_0$ after the
translation and hence becomes $\Omega_0$ after reversal and the
final application of $w_0$.

This fixes the compatible reversed type once and for all.
For $\de\in B_\mu^v(\lambda)$, the element $\lambda+\mu^*$ belongs
to the coroot lattice, so $t_\lambda$ and $t_{-\mu^*}$ induce the
same permutation of affine types. In word notation one must retain
the length-zero factor: merely reversing the simple reflections
would omit this relabeling.

\begin{prop}\label{prop:transported-duality}
\label{prop:reversal-criterion}
Literal reversal gives a bijection
\begin{equation}\label{eq:literal-duality}
B_\mu^{w_{0,\mu^*}}(\lambda)\longrightarrow B_{\mu^*}^{w_{0,\mu}}(-\lambda),
\qquad \de\longmapsto\de^*.
\end{equation}
It is the duality transported through the skeleton bijections, and
\[
(\de^*)^*=\de,\qquad
(e_\alpha\de)^*=f_\alpha(\de^*),\qquad
(f_\alpha\de)^*=e_\alpha(\de^*).
\]
Undefined operators are interpreted using $\mathbf0^*=\mathbf0$.
\end{prop}
\begin{proof}
A fold is positive precisely when its alcove lies on the negative
side of its wall. Translation and reversal preserve this property.
Thus $\de^*$ is positively folded.

Write the first and last alcoves as $u_0\Omega_0$ and
$t_{-\lambda}u\Omega_0$. Let $c_+,c_-$ be the numbers of positive
and negative crossings, and let $f$ be the number of folds.
For an alcove $A$, put
$\Psi(A)=-\sum_{\beta\in\Phi^+}
\lfloor\langle\beta,x\rangle\rfloor$, where $x\in A$;
this is independent of $x$. A positive crossing increases
$\Psi$ by one, a negative crossing decreases it by one,
and a fold leaves it unchanged. Since
$\Psi(u_0\Omega_0)=\ell(u_0)$ and
$\Psi(t_{-\lambda}u\Omega_0)=\langle2\rho,\lambda\rangle+\ell(u)$,
\[
c_+-c_-=\langle2\rho,\lambda\rangle+\ell(u)-\ell(u_0).
\]
Reversal interchanges positive and negative crossings and
preserves $f$. Thus
$\dim\de=\ell(u_0)+c_++f$ and
$\dim\de^*=\ell(u)+c_-+f$, giving
\begin{equation}\label{eq:reversal-dimension}
\begin{aligned}
\dim\de^*-\dim\de
&=\ell(u)-\ell(u_0)+c_--c_+\\
&=-\langle2\rho,\lambda\rangle.
\end{aligned}
\end{equation}
Consequently
\[
\dim\de^*=\langle\rho,\mu^*-\lambda\rangle+\ell(w_{0,\mu^*}),
\]
since $\langle\rho,\mu^*\rangle=\langle\rho,\mu\rangle$.
If the initial alcove of $\de^*$ has stabilizer factor $v'$
in the dual model, Lemma~\ref{lem:walk-dimension-bound} applied
to that model implies
\[
\ell(w_{0,\mu^*})\leq\ell(v')\leq\ell(w_{0,\mu})=\ell(w_{0,\mu^*}).
\]
Hence $v'=w_{0,\mu}$, and the reversed walk is maximal-dimensional.

For the chosen dual type, selecting skeleton faces commutes with
translation and reversal:
$\operatorname{sk}(\de^*)=(\operatorname{sk}(\de))^*$.
The gallery duality identities of \cite[Section 6]{GL05} and
uniqueness of the maximal lift now prove the operator identities.
Reversing twice gives the original walk, so the map is a bijection.
\end{proof}

\section{A Levi branching rule}
\label{sec:levi-branching}

We now fix a standard parabolic subgroup $P=P_J=MN$.
Our goal is to stratify the parabolic MV intersection
$K_Pt^\lambda K/K\cap Gr_{G,\mu}$ and use its maximal-dimensional
strata to prove the Levi branching rule.
We follow the geometric strategy of \cite{GL05}:
describe retraction fibers in the gallery variety, intersect
them with the open locus of minimal galleries, and count the
strata of maximal dimension. The additional step is
to show that the maximal-dimensional galleries lie in the
closed $M$-antidominant chamber. We prove this containment by an excursion reflection.

Throughout, $\mu$ is $G$-dominant and $\lambda$ is $M$-dominant.
We may assume $\mu-\lambda\in\mathbb Z\Phi^\vee$.
Otherwise the affine Grassmannian orbit intersections are
empty, galleries of the prescribed type cannot have weight
$\lambda$, and the branching multiplicity is zero.
Put
\[
\Phi_N=\Phi^+\setminus\Phi_M^+,
\qquad
\overline{\calC_M^-}=\bigcap_{j\in J}H_{\alpha_j,\leq0}.
\]
Recall that a gallery of weight $\nu$ ends at $-\nu$, and that
the center is $\De_{I_P}=t_\zeta\Omega_0=t^{-\zeta}\Omega_0$
for an $M$-central, sufficiently $N$-dominant $\zeta$.
Thus its $M$-coordinates lie in the base $M$-alcove, while its
$N$-coordinates are far positive. We write
$\dim_P=\dim_{\De_{I_P}}$ and
\[
\Ga_P^+(\nu,\ga_\mu)
=\Ga_P^+(\ga_\mu)\cap\Ga(\nu,\ga_\mu).
\]

\subsection{Parabolic retraction strata}

Let $X=X_2$ and let $g:X\to Gr_{G,\leq\mu}$ be the endpoint
morphism of Section~4. Its restriction to
$X^\circ=g^{-1}(Gr_{G,\mu})$ is an isomorphism; points of
$X^\circ$ are the minimal galleries of the prescribed type.
For a combinatorial gallery $\ga$, define
\[
\begin{aligned}
\Sigma_P(\ga)
&=\{\widetilde\ga\in X:
       \rho_{I_P,\calA}(\widetilde\ga)=\ga\},\\
S_P(\ga)&=g\bigl(\Sigma_P(\ga)\cap X^\circ\bigr).
\end{aligned}
\]
The retraction is applied to every face, so the incidence
conditions are included in this definition. We first construct
the locally closed structures on these fibers.

\begin{lem}[An attracting action for the parabolic center]
\label{lem:parabolic-action}
Choose an integral cocharacter $\eta_1$ such that
\[
\langle\alpha,\eta_1\rangle=0\quad(\alpha\in\Phi_M),
\qquad
\langle\beta,\eta_1\rangle>0\quad(\beta\in\Phi_N).
\]
Choose a rational point $x_0$ in the base alcove and a positive
integer $d$ with $\eta_2=dx_0\in X_*(T)$. For a sufficiently
large positive multiple $R$ of $d$, the action
\[
a_R(c)=\bigl(R\eta_1+\eta_2\bigr)(c)\,r(c^d),
\qquad r(c):t\longmapsto ct,
\]
has isolated fixed points indexed by the combinatorial
galleries, and
\begin{equation}\label{eq:two-stage-parabolic-stratum}
\Sigma_P(\ga)=\{x\in X:\lim_{c\to0}a_R(c)x=\ga\}.
\end{equation}
\end{lem}

\begin{proof}
Work over an algebraic closure $\bk$ of $k$. Put
$x_R=(R/d)\eta_1+x_0$ and $b=t_{(R/d)\eta_1}\Omega_0$.
The faces of galleries in $X(\bk)$ form a bounded subset $\calY$.
By \cite[Lemma~5.3]{Hai24}, for sufficiently large
$R\in d\mathbb Z_{>0}$ we have
\[
\rho_{b,\calA}|_{\calY}=\rho_{I_P,\calA}|_{\calY}.
\]

The Iwahori subgroup $I_b$ has root factors
$u_\beta(t^{q_\beta}\bk[[t]])$, where
$q_\beta=\lceil-\beta(x_R)\rceil$, and torus factor $T(\bk[[t]])$.
For $z(t)=\sum_{n\geq q_\beta}z_nt^n$ and
$h_T(t)\in T(\bk[[t]])$, conjugation gives
\[
\begin{aligned}
a_R(c)u_\beta(z(t))a_R(c)^{-1}
&=u_\beta\!\left(\sum_{n\geq q_\beta}
 z_nc^{d(\beta(x_R)+n)}t^n\right)\longrightarrow1,\\
a_R(c)h_T(t)a_R(c)^{-1}
&=h_T(c^dt)\longrightarrow h_T(0).
\end{aligned}
\]
All exponents in the first line are positive integers because
$x_R\in b$. Thus, if $h\in I_b(\bk)$ has torus factor $h_T(t)$,
its conjugate tends to $h_T(0)$. For an apartment face $F$,
$\rho_{b,\calA}^{-1}(F)=I_b(\bk)F$, and both $a_R$ and $T(\bk)$
fix $F$. Hence, whenever $hF$ is a face of a gallery in $X(\bk)$,
\[
\lim_{c\to0}a_R(c)(hF)=h_T(0)F
=F=\rho_{I_P,\calA}(hF).
\]
Applying this to every face proves
\eqref{eq:two-stage-parabolic-stratum}. It also identifies the
fixed points with the finite set $\Ga(\ga_\mu)$.

\end{proof}

\begin{prop}
\label{prop:parabolic-strata}
For every combinatorial gallery $\ga$, the fiber
$\Sigma_P(\ga)$ is an affine space with one coordinate for
each $\De_{I_P}$-load-bearing wall. Moreover,
\[
S_P(\ga)\ne\varnothing
\quad\Longleftrightarrow\quad
\ga\in\Ga_P^+(\ga_\mu).
\]
For such a gallery, $S_P(\ga)$ is a geometrically irreducible
locally closed subvariety and
\begin{equation}\label{eq:parabolic-stratum-dimension}
\dim S_P(\ga)=\dim_P(\ga).
\end{equation}
\end{prop}

\begin{proof}
Use the representatives $(v_0,\ldots,v_s)$ of $\ga$ and
$g_i=v_0\cdots v_{i-1}$ from Section~3, with $g_0=1$.
In the root-subgroup charts of \cite[Section 8]{GL05},
the $i$th factor has the form
\[
h_i=
\left(
\prod_{(\beta,n)}
u_\beta(z_{i,\beta,n}t^n)
\right)\dot v_i.
\]
The product is taken over all pairs $(\beta,n)\in\Phi\times\mathbb Z$
such that $(\beta+n)|_{F_{m_i}}=0$ and
$(\beta+n)(v_iF_{n_i})<0$.
The products are taken in fixed root orders.

The action calculated in Lemma~\ref{lem:parabolic-action}
is diagonal in these coordinates:
\[
a_R(c):
z_{i,\beta,n}
\longmapsto
c^{\,d(\beta+n)(g_i^{-1}x_R)}z_{i,\beta,n}.
\]
Indeed, the coordinate indexed by $\beta+n$ in the $i$th
factor corresponds to the affine root $g_i(\beta+n)$ in
the apartment. None of these exponents is zero.
The chart is $a_R$-stable and contains $\ga$, so it contains
every point whose $a_R$-limit is $\ga$. Therefore
\[
[h_0,\ldots,h_s]\in\Sigma_P(\ga)
\quad\Longleftrightarrow\quad
z_{i,\beta,n}=0
\quad\text{whenever}\quad
(\beta+n)(g_i^{-1}x_R)<0.
\]
All the remaining coordinates are arbitrary.

Since
$ 
\Ga_i'=g_iF_{m_i},
 \Ga_i=g_iv_iF_{n_i},
$ 
a remaining coordinate corresponds precisely to a wall
$ 
g_iH_{\beta,-n}$ 
containing $\Ga_i'$ and separating $\Ga_i$ from
$\De_{I_P}$. Consequently,
\[
\begin{aligned}
\Sigma_P(\ga)
&\cong
\prod_{i=0}^s
\prod_{(\beta,n)}
\mathbb A_k^1,\\
\dim\Sigma_P(\ga)
&=
\sum_{i=0}^s
\#\left\{
H\in\calH:
\Ga_i'\subset H,\ 
H\text{ separates }\Ga_i\text{ from }\De_{I_P}
\right\}\\
&=\dim_P(\ga).
\end{aligned}
\]
Here the $i$th factor of the first product is taken over the pairs $(\beta,n)$
occurring in the expression for $h_i$ that also satisfy
$(\beta+n)(g_i^{-1}x_R)>0$.

We next impose the condition that the gallery be minimal.
For $i\geq1$, cancelling the preceding factors gives
\[
\begin{aligned}
&(h_0\cdots h_{i-1})^{-1}
\bigl(
\widetilde\Ga_{i-1},
\widetilde\Ga_i',
\widetilde\Ga_i
\bigr)\\
&\hspace{3em}
=
\bigl(
F_{n_{i-1}},
F_{m_i},
h_iF_{n_i}
\bigr).
\end{aligned}
\]
By \cite[Section 9, Lemma 11 and Proposition 7]{GL05},
this local gallery is minimal precisely when
$
h_iL^+P_{n_i}
\in
L^+P_{n_{i-1}}\dot\tau_iL^+P_{n_i}/L^+P_{n_i}.
$
The two large faces of the local gallery have
opposite germs in the residue of $F_{m_i}$. Their
stabilizers are opposite parabolics in its reductive quotient,
and
$
L^+P_{n_{i-1}}\dot\tau_iL^+P_{n_i}/L^+P_{n_i}
=
I\dot\tau_iL^+P_{n_i}/L^+P_{n_i}.
$

For the fixed minimal gallery type, local minimality at every
$i\geq1$ is equivalent to minimality of the whole gallery;
see \cite[Remark 8]{GL05}. Thus
\[
\Sigma_P(\ga)\cap X^\circ
=
\left\{
[h_0,\ldots,h_s]\in\Sigma_P(\ga):
h_iL^+P_{n_i}\in
I\dot\tau_iL^+P_{n_i}/L^+P_{n_i}
\text{ for }1\leq i\leq s
\right\}.
\]

In the $i$th quotient, the coordinates of $\Sigma_P(\ga)$
parametrize the orbit
\[
\bigl(\dot g_i^{-1}I_P\dot g_i\cap L^+P_{m_i}\bigr)
\dot v_iL^+P_{n_i}/L^+P_{n_i}.
\]
The image of
$\dot g_i^{-1}I_P\dot g_i\cap L^+P_{m_i}$
in the reductive quotient is the Borel determined by the
projection of $g_i^{-1}\De_{I_P}$ to this residue.
Therefore the local minimality condition has a solution
if and only if
\[
\begin{aligned}
&\left(
\bigl(\dot g_i^{-1}I_P\dot g_i\cap L^+P_{m_i}\bigr)
\dot v_iL^+P_{n_i}/L^+P_{n_i}
\right)\\
&\hspace{3em}\cap
\left(
I\dot\tau_iL^+P_{n_i}/L^+P_{n_i}
\right)
\ne\varnothing.
\end{aligned}
\]
Equivalently,
\[
\dot\tau_iL^+P_{n_i}
\in
\overline{
\bigl(\dot g_i^{-1}I_P\dot g_i\cap L^+P_{m_i}\bigr)
\dot v_iL^+P_{n_i}/L^+P_{n_i}
}.
\]
To verify this equivalence, the second orbit is open and
\[
\lim_{c\to0}a_2(c)y=\dot\tau_iL^+P_{n_i}
\qquad
\left(
y\in I\dot\tau_iL^+P_{n_i}/L^+P_{n_i}
\right).
\]
The closure of the first orbit is $a_2$-stable.
Thus a nonempty intersection forces the displayed
closure inclusion. Conversely, that inclusion implies
that the open second orbit meets the first orbit,
which is dense in its closure.

The Bruhat closure relation in
$L^+P_{m_i}/L^+P_{n_i}$ identifies this inclusion with
the existence of the successive face reflections
\[
E^{(0)}=g_i\tau_iF_{n_i},
\qquad
E^{(j)}=s_{\beta_j,a_j}E^{(j-1)},
\qquad
E^{(m)}=g_iv_iF_{n_i},
\]
where
$
\Ga_i'\subset H_{\beta_j,a_j}$,
$
H_{\beta_j,a_j}
\text{ separates }E^{(j)}\text{ from }\De_{I_P}.
$
These are exactly the conditions of
Definition~\ref{def:positive-folding}.
The local conditions involve separate coordinate blocks,
and there is no minimality restriction on the first factor.
We conclude that
\[
\Sigma_P(\ga)\cap X^\circ\ne\varnothing
\quad\Longleftrightarrow\quad
\ga\in\Ga_P^+(\ga_\mu).
\]

When nonempty, $\Sigma_P(\ga)\cap X^\circ$ is an open
subvariety of the affine space $\Sigma_P(\ga)$.
It is therefore geometrically irreducible, with
\[
\dim\bigl(\Sigma_P(\ga)\cap X^\circ\bigr)
=
\dim\Sigma_P(\ga)
=
\dim_P(\ga).
\]
Finally, the isomorphism $g|_{X^\circ}$ restricts to
$
\Sigma_P(\ga)\cap X^\circ
\xrightarrow{\sim}S_P(\ga).
$
Hence $S_P(\ga)$ is locally closed and geometrically
irreducible, and
$
\dim S_P(\ga)=\dim_P(\ga).
$

\end{proof}

\subsection{From strata to branching multiplicities}

Write
\[
Z_{I_P}(\nu)=I_Pt^\nu K/K\cap Gr_{G,\mu},
\qquad Z_P(\lambda)=K_Pt^\lambda K/K\cap Gr_{G,\mu}.
\]
The preceding proposition identifies the dimensions and nonempty strata needed for the component count.

\begin{prop}[Retraction decomposition]
\label{prop:parabolic-retraction-decomposition}
There are finite partitions into locally closed subvarieties
\begin{align}
Z_{I_P}(\nu)
&=\bigsqcup_{\ga\in\Ga_P^+(\nu,\ga_\mu)}S_P(\ga),
\label{eq:parabolic-paving}\\
Z_P(\lambda)&=\bigsqcup_{\nu\in W_M\lambda}Z_{I_P}(\nu).
\label{eq:KP-decomposition}
\end{align}
The second union is over distinct elements of $W_M\lambda$;
its index set can also be written $W_M/W_{M,\lambda}$.
\end{prop}

\begin{proof}
Retraction commutes with the endpoint map. Its vertex fiber
over $-\nu$ is $I_Pt^\nu K/K$, by
\cite[Section 5.1]{Hai24} and \cite[Section 6]{HKM}.
Use $g|_{X^\circ}$ and
Proposition~\ref{prop:parabolic-strata} to obtain the first
partition. For the second partition, the Iwahori orbit
decomposition for $M$ gives
\[
\begin{aligned}
K_Pt^\lambda K/K
&=LNL^+Mt^\lambda K/K\\
&=\bigsqcup_{w\in W_M/W_{M,\lambda}}
LNI_Mt^{w\lambda}K/K
=\bigsqcup_{\nu\in W_M\lambda}I_Pt^\nu K/K;
\end{aligned}
\]
see \cite[Lemma~6.8]{HKM}. Intersecting with $Gr_{G,\mu}$
gives \eqref{eq:KP-decomposition}.
\end{proof}

Set
\[
\begin{aligned}
\mathscr G_P(\lambda)
&=\bigsqcup_{\nu\in W_M\lambda}\Ga_P^+(\nu,\ga_\mu),\\
\mathscr G_P^{\max}(\lambda)
&=\{\ga\in\mathscr G_P(\lambda):\dim_P(\ga)=\langle\rho,\lambda+\mu\rangle\}.
\end{aligned}
\]
All component counts below are geometric component counts.
If $\langle\rho,\lambda+\mu\rangle<0$, the dimension bound of \cite[Proposition 3.2]{Hai24}
makes $Z_P(\lambda)$ empty, and the stratum calculation makes
$\mathscr G_P(\lambda)$ empty as well. We interpret the
component set below as empty in this case.

\begin{prop}
\label{prop:parabolic-component-count}
Every $\ga\in\mathscr G_P(\lambda)$ has $\dim_P(\ga)\leq \langle\rho,\lambda+\mu\rangle$.
Closure in $Z_P(\lambda)$ gives a bijection
\[
\mathscr G_P^{\max}(\lambda)
\xrightarrow{\sim}\operatorname{Irred}_{\langle\rho,\lambda+\mu\rangle} Z_P(\lambda),
\qquad
\ga\longmapsto\overline{S_P(\ga)}.
\]
In particular,
\begin{equation}\label{eq:preliminary-branching}
[V_\mu^{\widehat G}:V_\lambda^{\widehat M}]
=\#\mathscr G_P^{\max}(\lambda).
\end{equation}
\end{prop}

\begin{proof}
By \cite[Proposition 3.2]{Hai24},
\[
\dim Z_P(\lambda)\leq \langle\rho,\lambda+\mu\rangle,
\qquad
\#\operatorname{Irred}_{\langle\rho,\lambda+\mu\rangle} Z_P(\lambda)
=[V_\mu^{\widehat G}:V_\lambda^{\widehat M}].
\]
Propositions~\ref{prop:parabolic-strata} and
\ref{prop:parabolic-retraction-decomposition} partition this
variety into nonempty irreducible locally closed strata of
dimensions $\dim_P(\ga)$. This gives the bound.
A $\langle\rho,\lambda+\mu\rangle$-dimensional stratum has a $\langle\rho,\lambda+\mu\rangle$-dimensional irreducible
closure, which is a component. Conversely, a component of
dimension $\langle\rho,\lambda+\mu\rangle$ must contain a dense stratum, since a finite
union of smaller-dimensional subsets cannot cover it. This proves the bijection and the formula.
\end{proof}

\subsection{Excursion reflection and the chamber condition}

The following auxiliary reflection is the parabolic version
of the dimension-increasing excursion operation. It is
not a raising or lowering operator: there is no translation of the tail.

\begin{dfn}
\label{def:parabolic-excursion}
Let $\alpha=\alpha_j$ for $j\in J$, and write
\[
\ga:\Ga_0'\subset\Ga_0\supset\cdots
          \subset\Ga_s\supset\Ga_{s+1}'.
\]
If $\ga\subset H_{\alpha,\leq0}$, leave
$e_{\alpha,P}(\ga)$ undefined. Otherwise let $k$ be its
first departure into $H_{\alpha,>0}$, and let $q>k$ be its
first subsequent small face in $H_{\alpha,0}$.
If there is no such face, set $q=s+2$. Reflect precisely
the large and small faces of indices $k\leq i<q$ by
$s_\alpha$, leaving the remaining faces unchanged.
The result is $e_{\alpha,P}(\ga)$.
\end{dfn}

The departure and returning faces lie in $H_{\alpha,0}$
and are fixed by the reflection. An immediate return is
allowed; if the excursion never returns, the terminal
vertex is reflected as well.

\begin{lem}
\label{lem:parabolic-excursion-endpoint}
If $\ga\in\mathscr G_P(\lambda)$ and $e_{\alpha,P}(\ga)$
is defined, then it has the same type and belongs to
$\mathscr G_P(\lambda)$. Its endpoint is unchanged for
a returning excursion and is $-s_\alpha\nu$ instead of
$-\nu$ for a nonreturning excursion. Moreover,
\begin{equation}\label{eq:parabolic-excursion-check}
\dim_P(e_{\alpha,P}(\ga))=\dim_P(\ga)+1.
\end{equation}
\end{lem}

\begin{proof}
We use the wall comparison in \cite[Section~6]{GL05}.
Put $\widetilde\ga=e_{\alpha,P}(\ga)$ and $H=H_{\alpha,0}$,
and use the indices $k,q$ of Definition~\ref{def:parabolic-excursion}.
For $b=t_\zeta\Omega_0$, the equality $s_\alpha\zeta=\zeta$
shows that $b$ and $s_\alpha b$ are adjacent across $H$.
Thus they lie on the same side of every wall other than $H$.
On the reflected portion, $H'\mapsto s_\alpha H'$ therefore
preserves load-bearing contributions and positive folding signs
for $H'\ne H$.

It remains to check $H$ at the two joins. At departure,
\[
\Ga_k'\subset H,\qquad
\Ga_k\subset H_{\alpha,>0},\qquad
\widetilde\Ga_k=s_\alpha\Ga_k\subset H_{\alpha,<0}.
\]
Hence $H$ becomes load-bearing at index $k$. It contains no
small face strictly inside the excursion, and the outgoing pair
at the return is unchanged. This proves
$\dim_P(\widetilde\ga)=\dim_P(\ga)+1$.

For positive folding, at $k>0$ append
$\Ga_k\xrightarrow{s_\alpha}s_\alpha\Ga_k$ to the original
folding sequence. At a return $q\leq s$, let $E$ be the original
minimal continuation of $\Ga_{q-1}$ across $\Ga_q'$.
Its germ is opposite to that of $\Ga_{q-1}$, so
$E\subset H_{\alpha,<0}$; prepend
$s_\alpha E\xrightarrow{s_\alpha}E$ to the original sequence.
Both added reflections are positive because their images lie in
$H_{\alpha,<0}$ and $b\subset H_{\alpha,>0}$.
There is no folding condition at an endpoint.
Thus $\widetilde\ga$ is positively folded.

The reflection fixes the joining faces and preserves face types. The endpoint stays $-\nu$ if $q\leq s+1$ and
becomes $-s_\alpha\nu$ otherwise; both belong to $-W_M\lambda$.
\end{proof}

\begin{prop}
\label{prop:parabolic-chamber}
Every $\ga\in\mathscr G_P^{\max}(\lambda)$ lies in
$\overline{\calC_M^-}$ and ends at $-\lambda$.
\end{prop}

\begin{proof}
If $\ga$ leaves $H_{\alpha,\leq0}$ for a simple root of
$M$, Lemma~\ref{lem:parabolic-excursion-endpoint} gives a
gallery in the same endpoint orbit with dimension $\langle\rho,\lambda+\mu\rangle+1$,
contradicting Proposition~\ref{prop:parabolic-component-count}.
Thus all the defining half-space conditions hold.
If the endpoint is $-\nu$, then $\nu\in W_M\lambda$
and chamber containment makes $\nu$ $M$-dominant.
The unique dominant representative of this orbit is
$\lambda$, so $\nu=\lambda$.
\end{proof}

\subsection{Comparison with LS galleries}

\begin{lem}
\label{lem:chamber-orientations}
For a gallery $\ga$ contained in $\overline{\calC_M^-}$,
\begin{align}
\ga\in\Ga_P^+(\ga_\mu)
&\quad\Longleftrightarrow\quad\ga\in\Ga^+(\ga_\mu),
\label{eq:parabolic-orientation-check}\\
\dim_P(\ga)&=\dim(\ga).
\label{eq:parabolic-dimension-check}
\end{align}
The two retraction cells in $X$ are also equal:
$\Sigma_P(\ga)=C_{X_2}(\ga)$.
\end{lem}

\begin{proof}
Consider a wall $H_{\beta,n}$ containing a small face of
$\ga$, and choose $\beta\in\Phi^+$. If $\beta\in\Phi_N$,
both centers are on its positive side, by sufficient
dominance. If $\beta\in\Phi_M^+$, chamber containment
gives $n\leq0$. The parabolic center satisfies
$0<\beta(b)<1$, whereas the Borel center has $\beta$
arbitrarily large. Both are therefore on the positive
side of this wall as well. This includes $n=0$.

Thus the centers induce exactly the same orientation in
the residue of every small face. This proves both the
positive-folding comparison and equality of the
load-bearing-wall counts. In the root-coordinate
description of Proposition~\ref{prop:parabolic-strata},
the same coordinates are free and the same coordinates
are zero at every step. Choosing the same root orders
therefore identifies the cells as the same subset of $X$.
\end{proof}

Put
\[
\mathscr L_M(\lambda)
=\{\ga\in\Ga_{LS}^+(\lambda,\ga_\mu):
                    \ga\subset\overline{\calC_M^-}\}.
\]

\Needspace{10\baselineskip}
\begin{thm}[Gallery branching rule]
\label{thm:levi-gallery-branching}
One has
\begin{equation}\label{eq:maximal-parabolic-LS}
\mathscr G_P^{\max}(\lambda)=\mathscr L_M(\lambda),
\end{equation}
and consequently
\[
[V_\mu^{\widehat G}:V_\lambda^{\widehat M}]
=\#\{\ga\in\Ga_{LS}^+(\lambda,\ga_\mu):
                    \ga\subset\overline{\calC_M^-}\}.
\]
\end{thm}

\begin{proof}
Proposition~\ref{prop:parabolic-chamber} puts every
maximal-dimensional parabolic gallery in the displayed
chamber and fixes its endpoint at $-\lambda$.
Lemma~\ref{lem:chamber-orientations} makes it Borel-positive
of dimension $\langle\rho,\lambda+\mu\rangle$, hence LS.
Conversely, the same lemma makes any gallery in
$\mathscr L_M(\lambda)$ parabolic-positive of dimension
$\langle\rho,\lambda+\mu\rangle$. This proves the equality of sets. The multiplicity
formula follows from
Proposition~\ref{prop:parabolic-component-count}.
\end{proof}

\begin{rmk}
This proof is the parabolic analogue of the cell and excursion arguments in \cite{GL05}. It gives a geometric proof of the gallery branching rule, independently of the path-model rule in \cite{Lit94}.
\end{rmk}

\subsection{Comparison with Haines' alcove walks}

Let $\mu$ be $G$-dominant and let $\lambda$ be $M$-dominant.
Fix the reduced expressions for $w(t_{\mu^*})_0$,
$w\in W^{\mu^*}$, chosen in Section~4.
Following \cite[Definition~5.8 and equation~(5.5)]{Hai24}, let
$\mathcal P_\mu^{\De_{I_P}}(\lambda)$ be the set of
$\De_{I_P}$-positively folded alcove-to-alcove walks of these
types, starting at $\Omega_0$ and ending at a vertex in
$-W_M\lambda$. The terminal vertex type is retained as in
\eqref{eq:extended-word}. Define
\[
\mathcal M_\mu^{\De_{I_P}}(\lambda)
=
\left\{
\delta'\in\mathcal P_\mu^{\De_{I_P}}(\lambda)
\;\middle|\;
\dim_P\delta'=\langle\rho,\lambda+\mu\rangle
\right\}.
\]
Here $\dim_P\delta'$ is the number of
$\De_{I_P}$-positive crossings plus the number of
$\De_{I_P}$-positive folds of $\delta'$.

\Needspace{10\baselineskip}
\begin{prop}[The parabolic skeleton bijection]
\label{prop:parabolic-walk-criterion}
After removing the initial minimal segment, the skeleton
map induces a bijection
\[
\mathcal M_\mu^{\De_{I_P}}(\lambda)
\xrightarrow{\sim}\mathscr L_M(\lambda).
\]
Its inverse is the Borel inverse skeleton map $L$ of
Section~\ref{subsec:haines-inverse}, with the initial
minimal segment prepended as in Section~4.
\end{prop}

\begin{proof}
Let $\ga\in\mathscr L_M(\lambda)$ and put
$\delta=L(\ga)=\operatorname{PR}_e^{-1}(\ga)\in B_\mu^e(\lambda)$,
using Proposition~\ref{prop:skeleton-bijection}.
If $w=t(\ga)$, then $\delta$ starts in $w\Omega_0$, ends at
$-\lambda$, and satisfies
$\dim\delta=\langle\rho,\lambda+\mu\rangle$.

We compare the two orientations on $\delta$.
The reduced block in \eqref{eq:type-blocks} lies in $W_{m_k}$,
so each of its simple reflections fixes $F_{m_k}$.
For $i_{k-1}<i\leq i_k$, the paragraph following Definition~\ref{def:skeleton}
therefore gives
\[
\Gamma_k'\subset\Delta_i'\subset H_{\beta_i,n_i},
\qquad \beta_i\in\Phi^+,
\]
also when the $i$th step is a fold.
If $\beta_i\in\Phi_M^+$, then
$\Gamma_k'\subset\overline{\calC_M^-}$ implies $n_i\leq0$.
Hence $\De_U$ and $\De_{I_P}$ both lie in
$H_{\beta_i,>n_i}$. For $\beta_i\in\Phi_N$, the same
inclusion follows from sufficient $N$-dominance.
The crossing and folding signs therefore agree at every step.
The initial-vertex contribution is $\ell(w)$ for both centers,
so $\delta$ is $\De_{I_P}$-positively folded and
$\dim_P\delta=\dim\delta$.

Prepend the fixed minimal segment from $\Omega_0$ to
$w\Omega_0$, as in Section~4, and call the resulting walk
$\delta'$. It has type $w(t_{\mu^*})_0$ and endpoint
$-\lambda$. The prepended segment has $\ell(w)$ positive
crossings, replacing the initial-vertex contribution of
$\delta$. Consequently,
\[
\dim_P\delta'=\dim_P\delta
=\langle\rho,\lambda+\mu\rangle,
\]
and $\delta'\in\mathcal M_\mu^{\De_{I_P}}(\lambda)$.

The map $\ga\mapsto\delta'$ is injective, since removing the
prepended segment and taking the skeleton recovers $\ga$.
By Theorem~\ref{thm:levi-gallery-branching} and
\cite[Theorem~5.14(2)]{Hai24},
\[
\#\mathscr L_M(\lambda)
=[V_\mu^{\widehat G}:V_\lambda^{\widehat M}]
=\#\mathcal M_\mu^{\De_{I_P}}(\lambda).
\]
Thus it is a bijection, with inverse given by the skeleton
after removal of the initial minimal segment.
\end{proof}

By definition, every alcove walk in $\mathcal M_\mu^{\De_{I_P}}(\lambda)$ has endpoint in the set $-W_M\la$. The previous proposition tells us that the endpoint has to be $-\la$.

\begin{cor}
    Every alcove walk in $\mathcal M_\mu^{\De_{I_P}}(\lambda)$ has endpoint $-\la$.
\end{cor}

\begin{cor}[Chamber containment for alcove walks]
\label{cor:parabolic-walk-chamber}
For every
$\gamma\in\mathscr L_M(\lambda)$, the walk
$\operatorname{PR}_{w_{0,\mu^*}}^{-1}(\gamma)$ is contained in
$\overline{\calC_M^-}$.
Consequently, if $\mu$ is regular, every walk in
$\mathcal M_\mu^{\De_{I_P}}(\lambda)$ is contained in this
closed chamber after removal of its initial minimal segment.
\end{cor}

\begin{proof}
Put $\delta=\operatorname{PR}_{w_{0,\mu^*}}^{-1}(\gamma)$ and
$w=t(\gamma)$. Write its alcoves as $\Delta_0,\ldots,\Delta_p$,
so that $\Delta_0=ww_{0,\mu^*}\Omega_0$ and
$\Gamma_0=w\Upsilon_0\subset\overline{\Delta_0}$.
Fix a simple root $\alpha$ of $M$.
Since $\gamma\subset H_{\alpha,\leq0}$ and $\Gamma_0'=\{0\}$,
Lemma~\ref{lem:skeleton-height} gives
$m_\alpha(\delta)=m_\alpha(\gamma)=0$.

If $\Gamma_0\subset H_{\alpha,<0}$,
then $\Delta_0\subset H_{\alpha,<0}$.
Otherwise $\Gamma_0\subset H_{\alpha,0}$, so
$\langle w^{-1}\alpha,\mu^*\rangle=0$.
Since $w\in W^{\mu^*}$ and $\alpha\in\Phi^+$,
\[
w^{-1}\alpha\in\Phi^+,\qquad
w_{0,\mu^*}^{-1}w^{-1}\alpha\in\Phi^-.
\]
Here $w_{0,\mu^*}$ sends the positive roots orthogonal to $\mu^*$
to negative roots. Thus $\Delta_0=ww_{0,\mu^*}\Omega_0\subset
H_{\alpha,<0}$ in this case as well.
Lemma~\ref{lem:maximal-wall}(2) now gives
$\Delta_i\subset H_{\alpha,<0}$ for every $i$.
Applying this to all simple roots of $M$ yields
\[
\bigcup_{i=0}^p\overline{\Delta_i}
\subset\bigcap_{j\in J}H_{\alpha_j,\leq0}
=\overline{\calC_M^-}.
\]

For $\delta'\in\mathcal M_\mu^{\De_{I_P}}(\lambda)$,
Proposition~\ref{prop:parabolic-walk-criterion} gives
its vertex-to-vertex remainder as
$\operatorname{PR}_e^{-1}(\gamma)$ for some
$\gamma\in\mathscr L_M(\lambda)$.
If $\mu$ is regular, then $W_{\mu^*}=\{e\}$ and $w_{0,\mu^*}=e$,
so the preceding containment applies to this remainder.
\end{proof}

\begin{rmk}

It is not necessarily true that every alcove walk in $\mathcal M_\mu^{\De_{I_P}}(\lambda)$ lies in $\overline{\calC_M^-}$ after removing the initial minimal walk. One can take a look at the $SL_3$ example in Figure~\ref{fig:A2-LS-lift}: If we take $M$ to be the Standard Levi corresponding to $\al_1$, we can see that even if $\ga$ lies entirely in the $M$-chamber $\overline{\calC_M^-}$, its alcove walk lift $\de$ does not.
\end{rmk}

\end{document}